\documentclass{elsarticle}
\usepackage{amssymb}
\usepackage{amsmath}
\usepackage{amsthm}
\usepackage{hyperref}
\hypersetup{
    linktocpage=true,
    colorlinks=true,
    linkcolor=blue,
    filecolor=magenta,      
    urlcolor=cyan,
    pdftitle={Multiplication operators on non-commutative spaces},
    bookmarks=true,
   citecolor=red    }
\newtheorem{theorem}{Theorem}[section]
\newtheorem{lemma}[theorem]{Lemma}
\newtheorem{proposition}[theorem]{Proposition}
\newtheorem{corollary}[theorem]{Corollary}
\newtheorem{remark}[theorem]{Remark}
\newtheorem{example}[theorem]{Example}

\newcommand{\dfs}[2]{d\left({#1}\right)({#2})}   

\newcommand{\svf}[1]{\mu_{{#1}}} 

\newcommand{\svft}[2]{\mu_{{#1}}({#2})}

\newcommand{\limm}[1]{\underset{#1}{\lim}}

\newcommand{\supu}[1]{\underset{#1}{\sup}\,}

\newcommand{\infu}[1]{\underset{#1}{\inf}\,}

\newcommand{\norm}[1]{\bigl\Vert{#1} \bigr\Vert} 

\newcommand{\normx}[2]{\bigl\Vert{#1} \bigr\Vert_{{#2}}}

\newcommand{\norminf}[1]{\bigl\Vert{#1} \bigr\Vert_{\mathcal{A}}} 

\usepackage{dsfont}
\newcommand{\id}{\mathds{1}} 

\newcommand{\summ}[2]{\underset{#1}{\overset{#2}{\sum}}}

\newcommand{\ip}[2]{\left\langle {#1},{#2} \right\rangle}

\newcommand{\dsumx}[2]{ \underset{#1}{\overset{#2}{\oplus}}}

\begin{document}

\title{Multiplication operators on non-commutative spaces}

\author[pdj]{P. de Jager\corref{cor1}}
\ead{28190459@nwu.ac.za}

\author[lel]{L. E. Labuschagne}
\ead{Louis.Labuschagne@nwu.ac.za}

\cortext[cor1]{Corresponding author}

\address[pdj]{DST-NRF CoE in Math. and Stat. Sci,\\ Unit for BMI,\\ Internal Box 209, School of Comp., Stat., $\&$ Math. Sci.\\
NWU, PVT. BAG X6001, 2520 Potchefstroom\\ South Africa}

\address[lel]{DST-NRF CoE in Math. and Stat. Sci,\\ Unit for BMI,\\ Internal Box 209, School of Comp., Stat., $\&$ Math. Sci.\\
NWU, PVT. BAG X6001, 2520 Potchefstroom\\ South Africa}


%\subjclass[2010]{}

\date{\today}

\begin{abstract}
Boundedness and compactness properties of multiplication operators on quantum (non-commutative) function spaces are investigated. For endomorphic multiplication operators these properties can be characterized in the setting of quantum symmetric spaces. For non-endomorphic multiplication operators these properties can be completely characterized in the setting of quantum $L^p$-spaces and a partial solution obtained in the more general setting of quantum Orlicz spaces. 
\end{abstract}

\maketitle

\section{Introduction}

In recent years there appears to be a renewed interest in the study of multiplication operators. Even in the commutative setting new results regarding multiplication operators on Orlicz spaces \cite{key-Chaw16,key-Chaw17}, Orlicz-Lorentz sequence spaces \cite{key-Bala13} and K\"othe sequence spaces \cite{key-Ramos17} have recently been obtained. In the non-commutative setting, multiplication operators have been studied on von Neumann algebras and their preduals, and between distinct Orlicz spaces \cite{key-Lab14}. In these articles sufficient conditions for the existence of multiplication operators between distinct Orlicz  spaces and necessary conditions for the compactness of multiplication operators between the respective spaces have been provided. 

It is also important to note the study of generalized duality, for which the underlying idea seems to be the identification of a particular space as the space of multipliers between two symmetric spaces. Numerous articles (\cite{key-Mali89}, \cite{key-Cala08} and \cite{key-Delgado10}, for example) have been written on this topic in the commutative setting and, more recently, some of these results have been generalized to the non-commutative setting (\cite{key-Han15} and \cite{key-Han16}). In particular, the space of multipliers between distinct non-commutative Calder\'{o}n-Lozanovski\u{i} spaces (generalizations of Orlicz-Lorentz spaces) is described, under the proviso that the (right continuous inverses of the) Orlicz functions satisfy certain inequalities. 

%%L
In this article we complement these results by characterizing the existence, boundedness and compactness of multiplication operators between distinct non-commutative Orlicz spaces, provided the Orlicz functions satisfy certain composition relations.  We choose to follow an approach focusing on the individual multiplier, rather than the identification of spaces of multipliers. (This decision is in part motivated by the fact that some questionable results pertaining to the non-commutative case have started appearing in the generalized duality literature - see \S 3.) The aforementioned conditions on the Orlicz functions engender a generalization of the setting of multiplication operators from an $L^p$-space into an $L^q$-space, where $p>q$. We will also characterize these properties for the case $p<q$ by using non-commutative analogues of the techniques employed in \cite{key-Takagi99}. Here our results clearly show that in this setting boundedness and compactness of a multiplication operator is dependent on the specific structure of the individual multiplier, and is not conditioned by membership of the multiplier to some a priori given space. It is therefore our belief that the ``generalized duality'' approach simply does not work in this setting. Regarding the endomorphic setting, we show that these properties can be characterized in the general setting of symmetric spaces. 
 
%%L
Throughout this paper we have confined ourselves to non-commutative spaces associated with semi-finite von Neumann algebras. The recent construction of Orlicz spaces for type III von Neumann algebras raises the intriguing possibility of ultimately extending the results herein to such spaces.

\section{Preliminaries}

Throughout this paper $\mathcal{A}$ will be used to denote a semi-finite von Neumann algebra equipped with a faithful normal semi-finite trace $\tau$. We will use $\id$ to denote the identity of $\mathcal{A}$. If $\mathcal{A}$ does not contain minimal projections, then it is called \emph{non-atomic}. A von Neumann algebra is called \emph{purely atomic} if it contains a set $\{p_\lambda\}_{\lambda \in \Lambda}$ of minimal projections such that $\sum p_\lambda =\id$, and this happens if and only if it is a product of Type 1 factors (see \cite[p.354]{key-Black06}). Furthermore, there exists a unique central projection $c\in \mathcal{A}$ such that $c\mathcal{A}$ is purely atomic and $c^\perp \mathcal{A}$ is non-atomic (this result follows from the corresponding result for $JBW$-algebras - see  \cite[Lemma 3.42]{key-Alfsen03}). The set of all $\tau$-measurable operators affiliated with $\mathcal{A}$ will be denoted $S(\mathcal{A},\tau)$. Let $x\in S(\mathcal{A},\tau)$ and let $|x|=\int_0^\infty \lambda de^{|x|}(\lambda)$ denote the spectral decomposition of $|x|$. We define the \emph{distribution function} of $|x|$ as \[\dfs{|x|}{s}:=\tau\left(e^{|x|}(s,\infty)\right) \qquad s\geq 0.\] The \emph{singular value function} of $x$, denoted $\svf{x}$, is defined to be the right continuous inverse of the  distribution function of $|x|$, namely \[\svft{x}{t}=\inf\{s\geq 0: \dfs{|x|}{s}\leq t\} \qquad t\geq 0. \]  This is the non-commutative analogue of the concept of a decreasing rearrangement of a measurable function. If $x,y \in S(\mathcal{A},\tau)$, then we will say that $x$ is \emph{submajorized} by $y$ and write $x \prec\prec y$ if 
\[\int_0^t \svft{x}{s}ds\leq \int_0^t \svft{y}{s}ds \qquad \text{for all $t> 0$.}\]
A linear subspace $E\subseteq S(\mathcal{A},\tau)$, equipped with a norm $\norm{\cdot}_E$, is called a \emph{symmetric space} if 
\begin{itemize}
\item $E$ is complete; 
\item $uxv\in E$ whenever $x\in E$ and $u,v\in\mathcal{A}$;
\item $\norm{uxv}_E\leq \norm{u}_{\mathcal{A}}\norm{v}_{\mathcal{A}}\norm{x}_E$ for all $x\in E,u,v\in \mathcal{A}$;
\item and $x\in E$ with $\norm{x}_E\leq \norm{y}_E$, whenever $y\in E$ and $x\in S(\mathcal{A},\tau)$ with $\svf{x}\leq \svf{y}$. 
\end{itemize}
It follows that $\norm{x}_E\leq \norm{y}_E$, whenever $E$ is a symmetric space and $x,y\in E$ with $|x|\leq |y|$. A symmetric space $E\subseteq S(\mathcal{A},\tau)$ is called \emph{strongly symmetric} if its norm has the additional property that $\norm{x}_E\leq \norm{y}_E$, whenever $x,y\in E$ satisfy $x \prec\prec y$.  If $E$ is a symmetric space and it follows from $x\in S(\mathcal{A},\tau)$, $y \in E$ and $x \prec\prec y$ that $x\in E$ and $\norm{x}_E \leq \norm{y}_E$, then $E$ is called a \emph{fully symmetric space}. Let $E \subseteq S(\mathcal{A},\tau)$ be a symmetric space. The \emph{carrier projection} $c_E$ of $E$ is defined to be the supremum of all projections in $\mathcal{A}$ that are also in $E$. If $c_E=\id$, then $E$ is continuously embedded in $S(\mathcal{A},\tau)$ equipped with the measure topology $\mathcal{T}_m$. We will therefore assume throughout this text that $c_E=\id$. Further details regarding $\tau$-measurable operators and symmetric spaces may be found in \cite{key-Terp1} and \cite{key-Dodds14}. We will focus on two particular examples of symmetric spaces, namely Orlicz spaces and $L^p$-spaces. 
 
A function $\varphi:[0,\infty)\rightarrow[0,\infty]$ is called an \emph{Orlicz (Young) function} if $\varphi$ is convex, $\varphi(0)=0$ and $\limm{t\rightarrow \infty} \varphi(t)=\infty$. We assume further that $\varphi$ is neither identically zero nor identically infinite on $(0,\infty)$ and that $\varphi$ is left continuous. Let 
\[a_{\varphi}:=\text{inf}\{t>0:\varphi(t)>0\} \qquad \text{and} \qquad b_{\varphi}:=\text{sup}\{t>0:\varphi(t)<\infty\}.\] Each Orlicz function $\varphi$ induces a \emph{complementary Orlicz function} $\varphi^*$ which is defined by $\varphi^*(s)=\supu{t>0}\{st-\varphi(t)\}$. The right continuous inverse of an Orlicz function $\varphi$ is defined by  $\varphi^{-1}(t):=\inf\{s:\varphi(s)>t\}=\sup\{s:\varphi(s)\leq t\}$. In the following proposition we present a few relevant properties of Orlicz functions, their complementary functions and right continuous inverses. 

\begin{proposition}\cite{key-BS88}\label{Rp.276 BS}\label{L4.8.16 BS}
Let $\varphi$ be an Orlicz function, $\varphi^*$ its complementary function and $\varphi^{-1}$ its right continuous inverse. Then
\begin{enumerate}
\item $\varphi(\varphi^{-1}(t))\leq t \leq \varphi^{-1}(\varphi(t))$,  for all  $0\leq t < \infty$; 
\item $t\leq \varphi^{-1}(t).(\varphi^*)^{-1}(t)\leq 2t$, for all $0\leq t <\infty$.
\end{enumerate}
\end{proposition}

 Suppose $(\Omega,\Sigma,\mu)$ is a measure space. If $\varphi$ is an Orlicz function, we can define a modular $I_\varphi$ on $L^0 (\Omega,\Sigma,\mu)=L^0(\mu)$, the space of all (equivalence classes) of measurable functions on $\Omega$, by setting 
\[I_{\varphi}(f):=\int_{\Omega}\varphi(|f(t)|)d\mu.\]
The collection of all $f\in L^{0}(\mu)$ such that $I_{\varphi}(\lambda f)<\infty$ for some $\lambda>0$ is called an \emph{Orlicz space} and is denoted by $L^\varphi(\mu)$. Restricted to $L^{\varphi}(\mu)$, the functional $\norm{\cdot}_{L^{\varphi}(\mu)} \colon L^{0}(\mu)\rightarrow[0,\infty)$ defined by \[\norm{f}_{L^{\varphi}(\mu)}=\text{inf}\{\lambda^{-1}:I_{\varphi}(\lambda f)\leq1\}\] is a norm, called the \emph{Luxemburg-Nakano norm}. Detailed investigations of Orlicz spaces and their properties may be found in \cite{key-Kras61} and \cite{key-Rao91}. Having defined Orlicz spaces in the commutative setting, we can use singular value functions to define non-commutative analogues of these spaces in the following way. If $(\mathcal{A},\tau)$ is a semi-finite von Neumann algebra and $E(0,\infty)\subseteq L^{0}(0,\infty)$ is a (fully) symmetric space, then the collection 
\[\{x\in S(\mathcal{A},\tau):\svf{x}\in E(0,\infty)\}\]
will be denoted $E(\tau)$ and is a (fully) symmetric space, when equipped with the norm (\cite[Theorems 4.2 and 4.5]{key-dP89}, \cite[p.218,219]{key-Pag})
\[\norm{x}_{E(\tau)}=\norm{\svf{x}}_{E(0,\infty)} \qquad x\in E(\tau).\]
In particular, since  $L^\varphi(0,\infty)$ is a rearrangement invariant Banach function space with the Fatou property, by  
\cite[Theorem 4.8.9]{key-BS88}; it follows (see \cite[p.202]{key-Pag}) that $L^\varphi(0,\infty)$ is fully symmetric and therefore $L^\varphi(\tau)$ is fully symmetric. When dealing with a non-commutative Orlicz space $L^\varphi(\tau)$, we will often use $\norm{\cdot}_\varphi$ to denote its norm, unless we wish to highlight the distinction between this norm and the corresponding norm in the commutative setting. The following results contain information to be used in the sequel and also show that non-commutative Orlicz spaces can be equivalently defined using a more direct approach. An important consideration in this approach is the fact that if $\varphi$ is an Orlicz function and $x\in S(\mathcal{A},\tau)$, then $\varphi(|x|)$ may not exist as an element of $S(\mathcal{A},\tau)$, if $b_\varphi<\infty$, and therefore care is required.

\begin{lemma}\cite{key-Lab13}\label{L2.1 Lab13}
Let $\varphi$ be an Orlicz function and $x\in S(\mathcal{A},\tau)$ a $\tau$-measurable element for which $\varphi(|x|)$ is again $\tau$-measurable. Extend $\varphi$ to a function on $[0,\infty]$ by setting $\varphi(\infty)=\infty$. Then $\varphi(\svf{x})=\svf{\varphi(|x|)}$ and $\tau(\varphi(|x|))=\int_0^{\infty}\varphi(\svft{x}{t})dt$. In particular, if $b_\varphi=\infty$, then $\varphi(|x|)\in S(\mathcal{A},\tau)$ for all $x\in S(\mathcal{A},\tau)$.
\end{lemma}

\begin{lemma}\label{L1 20/09/17}
If $\varphi$ is an Orlicz function, then $\varphi^{-1}(|x|)\in S(\mathcal{A},\tau)$ whenever $x\in S(\mathcal{A},\tau)$.
\begin{proof} 
 Since $\limm{s\rightarrow \infty}\varphi(s)=\infty$, the set $\{s\geq 0: \varphi(s)>t\}$ is non-empty for each $t\geq 0$ and hence $\varphi^{-1}(t)=\inf \{s\geq 0: \varphi(s)>t\}$ is finite for each $t\geq 0$. Since $\varphi^{-1}$ is also increasing, this implies that $\varphi^{-1}$ is bounded on compact subsets of $[0,\infty)$. We therefore obtain $\varphi^{-1}(|x|)\in S(\mathcal{A},\tau)$ whenever $x\in S(\mathcal{A},\tau)$ (see \cite[Proposition 4.8]{key-Pag}). 
\end{proof}
\end{lemma}

\begin{proposition}\cite{key-Lab13} \label{P2.2 Lab13}
Let $\varphi$ be an Orlicz function and $x\in S(\mathcal{A},\tau)$. There exists some $\alpha>0$ such that $\int_0^\infty \varphi(\alpha\svft{x}{t})dt<\infty$ if and only if there exists some $\beta>0$ such that $\varphi(\beta|x|)\in S(\mathcal{A},\tau)$ and $\tau(\varphi(\beta|x|))<\infty$. Moreover \[\norm{\svf{x}}_{L^{\varphi}(0,\infty)}=\inf\{\lambda>0:\varphi(|x|/\lambda)\in S(\mathcal{A},\tau),\tau(\varphi(|x|/\lambda))\leq 1\}.\]
\end{proposition}

\begin{remark}\label{RP2.2 Lab13}
It is useful to note that in the proof of Proposition \ref{P2.2 Lab13} it is shown that if $x\in S(\mathcal{A},\tau)$ and $\alpha>0$ is such that $\int_0^{\infty}\varphi(\alpha \svft{x}{t})dt<\infty$, then for every $\epsilon>0$, $\varphi\left(\tfrac{\alpha}{1+\epsilon}x\right)\in \mathcal{A} \subseteq S(\mathcal{A},\tau)$ and by Lemma \ref{L2.1 Lab13}
\[\varphi\left(\tfrac{\alpha}{1+\epsilon}x\right)=\int_0^{\infty}\varphi\left(\tfrac{\alpha}{1+\epsilon}\svft{x}{t}\right)dt.\]
\end{remark}

We briefly mention K\"othe duality. Suppose $E\subseteq S(\mathcal{A},\tau)$ is a symmetric space. The collection 
\[E^{\times} :=\{x\in S(\mathcal{A},\tau):\tau(|yx|)<\infty\,\, \forall y\in E\}\]
is a symmetric space, called the \emph{K\"othe dual} of $E$, when equipped with the norm 
\[\norm{x}_{E^\times}:=\text{sup}\,\{\tau(|xy|):y\in E,\norm{y}_E\leq 1\}.\]It is known (see \cite{key-DDP2}) that 
\[E^{\times} =\{x\in S(\mathcal{A},\tau):yx\in L^1(\tau)\,\, \forall y\in E\}=\{x\in S(\mathcal{A},\tau):xy\in L^1(\tau)\,\, \forall y\in E\}\]
and if $E(0,\infty)$ is a (commutative) symmetric space, then 
\[(E(\tau))^\times=E^\times(\tau):=\{x\in S(\mathcal{A},\tau): \svf{x}\in E^\times(0,\infty)\},\]
where $E^\times(0,\infty):=\{f\in L^0(0,\infty):fg\in L^1(0,\infty)\, \forall g\in E(0,\infty)\}\}$. In the context of Orlicz spaces, we can identify the K\"othe dual as described in the following result.

\begin{proposition}\cite{key-Lab13}\label{P2.3 Lab13}
Let $\varphi$ be an Orlicz function and $\varphi^*$ its complementary function. Then $L^{\varphi^*}(\tau)$, equipped with the norm $\norm{\cdot}_{\varphi^*}^0$ defined for $x\in L^{\varphi^*}(\tau)$ by 
\begin{eqnarray*}
\norm{x}_{\varphi^*}^0&=&\sup\{\tau(|xy|):y\in L^\varphi(\tau),\norm{y}_{\varphi} \leq 1\} \\
&=&\infu{k>0}\left(\tfrac{1}{k}+\tfrac{1}{k}\int_0^\infty\varphi^*(k\svft{x}{t})dt\right),
\end{eqnarray*}
 is the K\"othe dual of $L^{\varphi}(\tau)$. Consequently \[|\tau(xy)|\leq \norm{x}_{\varphi^*}^0\norm{y}_{\varphi} \qquad \forall x\in L^{\varphi^*}(\tau), y\in L^\varphi(\tau).\]
\end{proposition}

\begin{remark}\label{RP2.3 Lab13}
If $E\subseteq S(\mathcal{A},\tau)$ is a symmetric space, then using the definition of $\norm{y}_{E^\times}$, it is easily verified that
\[\tau(|xy|)\leq \norm{x}_E \norm{y}_{E^\times},\]
whenever $x\in E$ and $y\in E^\times$. Since $|\tau(xy)|\leq \tau(|xy|)$, we obtain a sharper claim than the one made in Proposition \ref{P2.3 Lab13}.
\end{remark}

Next, we describe several growth conditions that will enable us to distinguish various classes of Orlicz spaces. The first such condition is the $\Delta_2$-condition. If there exists a $t_{0}>0$ and a $C>0$ such that $\varphi(2t)\leq C\varphi(t)<\infty$ for all $t$ such that $t_{0}\leq t<\infty$, then $\varphi$ is said to satisfy the \emph{$\Delta_{2}$-condition} for large $t$. If $t_{0}=0$, then $\varphi$ is said to satisfy the $\Delta_{2}$-condition globally and we write $\varphi\in\Delta_{2}$. The following details important consequences of the $\Delta_2$-condition.

\begin{proposition} \cite{key-BS88}\label{P invertible orlicz}
Suppose $\varphi$ is an Orlicz function. If  $\varphi\in \Delta_{2}$, then $\varphi$ is invertible  and for any $k>0$, there exists $m_{k}>0$ such that $\varphi(kt)\leq m_{k}\varphi(t)$, for all $t\geq 0$.
\end{proposition}

An Orlicz function $\varphi$ is said to satisfy the $\nabla'$-condition, if there exists a $t_{0}>0$ and a $c>0$ such that $\varphi(s)\varphi(t)\leq \varphi(cst)$ for all $s,t \geq t_0$. If $t_0=0$, then this condition is said to hold globally and we write $\varphi \in \nabla'$. We will be particularly interested in the following consequence of the $\nabla'$-condition.

\begin{lemma} \label{L2 26/09/2017}
Suppose $\varphi$ is an invertible Orlicz function. If $\varphi\in \nabla'$, then  
\[\varphi^{-1}(uv)\leq c \varphi^{-1}(u)\varphi^{-1}(v), \qquad \mbox{for all } u,v \geq 0, \]
where $c>0$ is such that $\varphi(s)\varphi(t)\leq \varphi(cst)$ for all $s,t \geq 0$.
\begin{proof}
Let $\epsilon>0$, $u\geq 0$ and $v \geq 0$ be given. Since each of $\varphi^{-1}(u)$ and $\varphi^{-1}(v)$ are finite, we may by the definition of $\varphi^{-1}$ select $r_1, r_2>0$ so that $$\varphi(r_1)>u, \quad \varphi(r_2)>v, \quad r_1\leq \varphi^{-1}(u)+\epsilon, \mbox{ and }r_2\leq \varphi^{-1}(v)+\epsilon.$$But since $\varphi(cr_1r_2)\geq \varphi(r_1)\varphi(r_2)>uv$, we must have that
$$\varphi^{-1}(uv)=\inf\{r>0: \varphi(r)>uv\}\leq cr_1r_2\leq c(\varphi^{-1}(u)+\epsilon)(\varphi^{-1}(v)+\epsilon).$$In view of the fact that $\epsilon$ was arbitrary, the claim follows.
\end{proof}
\end{lemma}

The non-commutative $L^p$-spaces can be defined as the collection of $\tau$-measurable operators whose singular value functions are $p$-integrable or equivalently as those $\tau$-measurable operators $x$ for which $\tau(|x|^p)<\infty$. Equipped with the norm 
\[\norm{x}_{L^p(\tau)}=\tau(|x|^p)^{1/p}=\int_0^\infty \left(\svft{x}{t}\right)^p dt, \qquad x\in L^p(\tau),\]
$L^p(\tau)$ is a symmetric space. Furthermore, we note that if $\varphi(t)=t^p$, for $t\geq 0$, then $\varphi$ is an Orlicz function, satisfying the $\Delta_{2}$- and $\nabla'$-conditions globally, and $L^\varphi(\tau)=L^p(\tau)$, with equality of norms. Unless confusion is possible, we will often denote the norm of an $L^p$-space using $\norm{\cdot}_p$.  If $1<p<\infty$, then we will use $p'$ to denote the conjugate index of $p$, i.e. $1/p+1/p'=1$. The following collects some of the relevant properties of $L^p$-spaces to be used in the sequel.

\begin{proposition}\cite{key-DDP2,key-dP89,key-Fack86} \label{T4.2 Fack86} \label{P3.4 DDP2}
Suppose $x,y\in S(\mathcal{A},\tau)$. Then
\begin{enumerate}
\item $\tau(xy)=\tau(yx)$, whenever $xy,yx\in L^1(\tau)$. If, in addition, $x,y\geq 0$, then $x^{1/2}yx^{1/2},y^{1/2}xy^{1/2}\in L^1(\tau)$ and \[\tau(xy)=\tau(x^{1/2}yx^{1/2})=\tau(y^{1/2}xy^{1/2});\]
\item $\normx{xy}{q}\leq \normx{x}{p}\normx{y}{r}$ whenever $p,q,r>0$ are such that $p^{-1}+r^{-1}=q^{-1}$; and
\item $\int_0^t f(\svft{xy}{s})ds\leq \int_0^t f(\svft{x}{s}\svft{y}{s})ds$ for any increasing function $f:\mathbb{R}^+\rightarrow \mathbb{R}$ such that $t\mapsto f(e^t)$ is convex.
\end{enumerate}
\end{proposition}

Suppose $E,F\subseteq S(\mathcal{A},\tau)$ are symmetric spaces and $w\in S(\mathcal{A},\tau)$. The left multiplication map  $E\rightarrow S(\mathcal{A},\tau) \colon x \mapsto wx$   will be denoted $M_w$. If $M_w$ maps $E$ into $F$, then $M_w$ will be called a \emph{multiplication operator} from $E$ into $F$. There are several natural questions regarding such multiplication maps. Firstly, what are the conditions on $w\in S(\mathcal{A},\tau)$ which characterize when $M_w$ maps $E$ into $F$? Furthermore, under what conditions will such multiplication operators be bounded or compact? Unsurprisingly, it is often the case that continuity properties of $M_w$ and conditions under which $M_w$ maps $E$ into $F$  are studied concurrently. In fact, $M_w$ is automatically continuous if it maps $E$ into $F$. (To see this observe that any $w\in S(\mathcal{A},\tau)$ induces a continuous (left) multiplication operator on $S(\mathcal{A},\tau)$. On combining this fact with the fact that each of $E$ and $F$ continuously embed into $S(\mathcal{A},\tau)$, it is now a simple exercise to show that $M_w$ must then have a closed graph as a map from $E$ to $F$.)

\section{Existence and boundedness of multiplication operators}

It is easily checked that  $w\in S(\mathcal{A},\tau)$ induces a bounded (left) multiplication operator if and only if $|w|$ induces a bounded (left) multiplication operator. Furthermore, if this is the case, then \[\norm{M_w}=\norm{M_{|w|}}.\] It therefore suffices to consider positive elements in our study of boundedness properties of multiplication operators. For the endomorphic setting the boundedness of multiplication operators has been characterized in the general setting of symmetric spaces (\cite[Proposition 5]{key-Han15}). For the non-endomorphic case, we will show that the boundedness of multiplication operators between different Orlicz spaces can be characterized if the Orlicz functions satisfy certain properties. These properties imply that the situation is a natural generalization of considering multiplication operators from  $L^p(\tau)$ into $L^q(\tau)$ if $p>q$. For $p<q$, a characterization will also be provided.

\subsection{Multiplication operators on symmetric spaces}

It is natural to consider if it is not possible to lift results from the commutative setting to the non-commutative setting. It is, in fact, claimed in \cite[Corollary 3.1]{key-Han16} that if $0<\alpha_0,\alpha_1<\infty$, $E, F\subseteq L^0(I)$ ($I=(0,1)$ or $I=(0,\infty)$) with $E$ an $\alpha_0$-convex symmetric quasi-Banach space and $F$ an $\alpha_1$-convex fully symmetric quasi-Banach space with the Fatou property, then $E(\tau)^{F(\tau)}$, the collection of all multipliers from $E(\tau)$ to $F(\tau)$, is given by 
\[E^F(\tau)=\{x\in S(\mathcal{A},\tau): \svf{x} \in E^F\},\]
where $E^F$ is the set of all multipliers from $E$ to $F$. Whilst there are other interesting and useful results in \cite{key-Han16}, the aforementioned result cannot be true without further restrictions on the semi-finite von Neumann algebra $\mathcal{A}$. On 
noting that symmetric Banach spaces are $1$-convex (see \cite{key-DDS14} for definitions and details) and $L^p$-spaces are fully symmetric spaces with the Fatou property, this can be seen from the following example. 

\begin{example}\label{E1 13/07/18}
Suppose $\mathcal{A}=\mathcal{B}(H)$, the set of all bounded operators on the Hilbert space $H$, and equip $\mathcal{A}$ with the canonical trace $tr$. If $1\leq p<q$, then $L^p(tr)$ is continuously embedded in $L^q(tr)$, by \cite[Proposition 4.5]{key-Diestel95}. It follows that the identity operator is a multiplier from $L^p(tr)$ into $L^q(tr)$. However, it follows from \cite[Theorem 1.4]{key-Takagi99} that the only multiplier from $L^p(I)$ to $L^q(I)$ is the function which is zero almost everywhere. It follows that $E^F(\tau)=\{0\}$. However as can be seen from Theorem \ref{T1b 13/10/17}, the set of multipliers yielding bounded operators $L^p(tr)$ to $L^q(tr)$, is actually quite large.
\end{example}

It is our suspicion that the discussion on p.286 (\cite{key-Fack86}) regarding the embedding of a general semi-finite von Neumann algebra $\mathcal{A}$ into the non-atomic von Neumann algebra $\mathcal{A}\bar{\otimes} L^\infty(0,\infty)$ has, at times, not been applied with sufficient care. This has led to an insufficient distinction between the atomic and non-atomic cases and possible mistakes in the literature. Regarding the proof of \cite[Corollary 3.1]{key-Han16}, it is true that $x$ and $x \otimes \id$ have the same generalized singular value functions. It need not, however, be the case that $x \otimes \id$ is a multiplier between two spaces corresponding to $\mathcal{A}\bar{\otimes}L^\infty(0,\infty)$ if $x$ is a multiplier between the matching spaces corresponding to $\mathcal{A}$. Example \ref{E1 13/07/18} above demonstrates that in the case $\mathcal{A}=\mathcal{B}(H)$, $\id$ is a multiplier from $L^p(tr)$ into $L^q(tr)$. However $\id \otimes \id$ is not a multiplier from $L^p(\mathcal{A}\bar{\otimes}L^\infty(0,\infty), tr \otimes m)$ into $L^q(\mathcal{A}\bar{\otimes}L^\infty(0,\infty), tr \otimes m)$.
To see this let $e$ be any minimal projection in $\mathcal{B}(H)$, $E_k=(\frac{1}{2^{k}},\frac{1}{2^{k-1}})$ for $k\in \mathbb{N}^+$ and $f_n=\summ{k=1}{n}\frac{\chi_{E_k}}{m(E_k)^{1/q}}$ for $n\in \mathbb{N}^+$. It is easily checked that $(f_n)_{n=1}^{\infty}$ is Cauchy in $L^p(0,\infty)$, but not in $L^q(0,\infty)$. Furthermore,  
\begin{eqnarray*}
e\otimes f_n &\in& (L^\infty\cap L^1)(\mathcal{A}\bar{\otimes}L^\infty(0,\infty), tr \otimes m) \\
& \subseteq& L^p(\mathcal{A}\bar{\otimes}L^\infty(0,\infty), tr \otimes m) \cap L^q(\mathcal{A}\bar{\otimes}L^\infty(0,\infty), tr \otimes m)
\end{eqnarray*}
and $\norm{(e \otimes f_n)-(e \otimes f_m)}_r=\norm{e \otimes (f_n-f_m)}_r=\norm{f_n-f_m}_r$ for every $n,m\in \mathbb{N}^+$ and $1\leq r <\infty$. It follows that $(e \otimes f_n)_{n=1}^\infty$ is Cauchy in $L^p(\mathcal{A}\bar{\otimes}L^\infty(0,\infty), tr \otimes m)$, but not in $L^q(\mathcal{A}\bar{\otimes}L^\infty(0,\infty), tr \otimes m)$. Since $(\id \otimes \id)(e\otimes f_n)=e \otimes f_n$ for each $n$, this shows that $\id \otimes \id$ is not a continuous multiplier from $L^p(\mathcal{A}\bar{\otimes}L^\infty(0,\infty), tr \otimes m)$ into $L^q(\mathcal{A}\bar{\otimes}L^\infty(0,\infty), tr \otimes m)$.

Cognizant of the above subtleties, we will not attempt to deal with the atomic case by means of a reduction to the non-atomic setting, but will rather follow a more direct approach. 

The next result shows that sufficient conditions for the existence and boundedness of multiplication operators on non-commutative spaces may however often be derived from the classical setting without imposing further restrictions on the semi-finite von Neumann algebra.

\begin{lemma}\label{P2 04/10/17}
Suppose $E_i(0,\infty)$ ($i=1,3$) are (commutative) symmetric spaces, $E_2(0,\infty)$ is a (commutative) strongly symmetric space and $(\mathcal{A},\tau)$ is a semi-finite von Neumann algebra. If there exists some $k>0$ such that
\begin{eqnarray}
\norm{fg}_{E_2(0,\infty)}\leq k\norm{f}_{E_1(0,\infty)}\norm{g}_{E_3(0,\infty)}, \label{e1 P2 04/10/17}
\end{eqnarray}
whenever $f\in E_1(0,\infty)$ and $g\in E_3(0,\infty)$, then 
\[\norm{xy}_{E_2(\tau)}\leq k\norm{x}_{E_1(\tau)}\norm{y}_{E_3(\tau)},\]
whenever $x\in E_1(\tau)$ and $y\in E_3(\tau)$.
\begin{proof}
Suppose $x\in E_1(\tau)$ and $y\in E_3(\tau)$. Then $\svf{xy}\prec\prec \svf{x}\svf{y}$, by \cite[Theorem 4]{key-Suk16}. Since $E_2(0,\infty)$ is a strongly symmetric space, this implies that $\norm{\svf{xy}}_{E_2(0,\infty)}\leq \norm{\svf{x}\svf{y}}_{E_2(0,\infty)}$. Using (\ref{e1 P2 04/10/17}) we therefore obtain
\begin{eqnarray*}
\norm{\svf{xy}}_{E_2(0,\infty)} \leq \norm{\svf{x}\svf{y}}_{E_2(0,\infty)} \leq k\norm{\svf{x}}_{E_1(0,\infty)} \norm{\svf{y}}_{E_3(0,\infty)}.
\end{eqnarray*} 
Since the norm of a trace-measurable operator is given by the norm of its singular value function, the result follows. 
\end{proof}
\end{lemma}

Regarding the endomorphic setting, we note that it is shown in \cite[Proposition 5]{key-Han15} that if $E\subseteq S(\mathcal{A},\tau)$ is a non-trivial symmetric space and $w\in S(\mathcal{A},\tau)^+$, then $M_w$ is a bounded multiplication operator from $E$ into itself if and only if $w\in \mathcal{A}$. Examination of the proof of this theorem shows that if this is the case, then $\norm{M_w}=\norm{w}_{\mathcal{A}}$.

\subsection{Multiplication operators on Orlicz spaces}

We note that it is shown in \cite[Theorem 3]{key-Han15}, that if $E$ is a symmetric space with the Fatou property, $\varphi, \varphi_1, \varphi_2$ are Orlicz functions with $b_\varphi=b_{\varphi_1}=b_{\varphi_2}$, and $k_1\varphi^{-1}\leq \varphi_1^{-1}(t)\varphi_2(t)\leq k_2\varphi^{-1}(t)$ for all $t\geq 0$ for some $k_1,k_2>0$, then the space of multipliers between the Calder\'{o}n-Lozanovski\u{i} spaces $E_{\varphi_1}(\tau)$ and $E_\varphi(\tau)$ is given by $E_{\varphi_2}(\tau)$. If $E=L^1$, then we obtain the corresponding result for Orlicz spaces. We extend this result in the setting of Orlicz spaces by showing that it holds true \textbf{without} the requirement that $b_\varphi=b_{\varphi_1}=b_{\varphi_2}$. We use the fact that the corresponding result holds true in the commutative setting (see \cite{key-Mali10}) to establish the first part of this result.

\begin{theorem}\label{T1 07/09/17}
Suppose $w\in S(\mathcal{A},\tau)^+$ and $\varphi_i$ $(i=1,2,3)$ are Orlicz functions.
\begin{enumerate}
\item If $\varphi^{-1}_1(t)\varphi_3^{-1}(t)\leq k\varphi_2^{-1}(t)$ for all $t\geq 0$, then $M_w$ is a bounded multiplication operator from $L^{\varphi_1}(\tau)$ into $L^{\varphi_2}(\tau)$ whenever $w\in L^{\varphi_3}(\tau)$. Furthermore, if this is the case, then $\norm{M_w}\leq 2k\norm{w}_{\varphi_3}$.
\item If there exists some $k>0$ such that $\varphi_2^{-1}(t)\leq k \varphi^{-1}_1(t)\varphi_3^{-1}(t)$ for all $t\geq 0$, then $M_w$ is a bounded multiplication operator from $L^{\varphi_1}(\tau)$ into $L^{\varphi_2}(\tau)$ only if $w\in L^{\varphi_3}(\tau)$. Furthermore, if this is the case, then $\norm{w}_{\varphi_3}\leq 4k \norm{M_w}$.
\end{enumerate}
\begin{proof}
(1) follows from \cite[Remark 2]{key-Mali10} and Lemma \ref{P2 04/10/17}. 

To prove (2), suppose that there exists some $k>0$ such that $\varphi_2^{-1}(t)\leq k \varphi^{-1}_1(t)\varphi_3^{-1}(t)$ for all $t\geq 0$ and $M_w$ is a bounded multiplication operator from $L^{\varphi_1}(\tau)$ into $L^{\varphi_2}(\tau)$. Let $\Gamma(\cdot):=\tau(w\cdot)$. We show that $\Gamma$ is a bounded linear functional on $L^{\varphi_3^*}(\tau)$. Let $\epsilon>0$ be given and suppose $x\in L^{\varphi_3^*}(\tau)^+$ with $\norm{x}_{\varphi_3^*}=1-\epsilon$. It follows by Proposition \ref{P2.2 Lab13} and Remark \ref{RP2.2 Lab13} that $\varphi_3^*(x)\in S(\mathcal{A},\tau)$ and $\tau(\varphi_3^*(x))\leq 1$. Define $x_1:=\varphi^{-1}_1\circ \varphi_3^*(x)$ and $x_2:=(\varphi_2^*)^{-1}\circ \varphi_3^*(x)$. Then $x_1,x_2 \in S(\mathcal{A},\tau)$, by Lemma \ref{L1 20/09/17}. Furthermore, $ \varphi_1 \circ \varphi_1^{-1} (t)\leq t$ for all $t\geq 0$, by Proposition \ref{Rp.276 BS}(1), and so $0\leq \varphi_1(x_1)=\varphi_1 \circ \varphi_1^{-1}(\varphi^*_3(x))\leq \varphi_3^*(x)$ using the Borel functional calculus. It follows that $\tau(\varphi_1(x_1))\leq \tau(\varphi_3^*(x))\leq 1$ and therefore,  $x_1 \in L^{\varphi_1}(\tau)$ with $\norm{x_1}_{\varphi_1}\leq 1$, by Proposition \ref{P2.2 Lab13}. Similarly, $x_2 \in L^{\varphi_2^*}(\tau)$ and $\norm{x_2}_{\varphi_2^*}\leq 1$. Since $M_w$ is a bounded multiplication operator from $L^{\varphi_1}(\tau)$ into $L^{\varphi_2}(\tau)$, we have that $wx_1\in L^{\varphi_2}(\tau)$. So $wx_1x_2\in L^1(\tau)$ by K\"othe duality. Applying this and Proposition \ref{P2.3 Lab13} (see also Remark \ref{RP2.3 Lab13}), we obtain
\begin{eqnarray}
&|\tau(wx_1x_2)|  \leq \tau(|wx_1x_2|)\leq \norm{wx_1}_{\varphi_2}\norm{x_2}_{\varphi_2^*}^0 \leq \norm{M_w}\norm{x_1}_{\varphi_1}\norm{x_2}_{\varphi_2^*}^0. & \label{e3 T1 25/08/17}
\end{eqnarray}
Furthermore, using \cite[Theorem 4.8.14]{key-BS88} and the fact that the norm of $y$ is equal to the norm of its singular value function in the corresponding commutative space, we have $\norm{y}_{\varphi_2^*}^0\leq 2\norm{y}_{\varphi_2^*}$ for any $y\in L^{\varphi_2^*}(\tau)$. Hence inequality \ref{e3 T1 25/08/17} becomes \begin{eqnarray} 
&|\tau(wx_1x_2)|  \leq 2\norm{M_w}\norm{x_1}_{\varphi_1}\norm{x_2}_{\varphi_2^*} \leq 2\norm{M_w}=\tfrac{2}{1-\epsilon}\norm{M_w}\norm{x}_{\varphi_3^*},& \label{e2 T1 25/08/17}
\end{eqnarray}
since $\norm{x_1}_{\varphi_1},\norm{x_2}_{\varphi_2^*}\leq 1$, as shown earlier and $\norm{x}_{\varphi_3^*}=1-\epsilon$. Since $x_1$ and $x_2$ were defined using $x$ and the Borel functional calculus, it is easily checked that $x_1x_2=x_2x_1\geq 0$. It follows that 
\[x_1x_2w=x_2x_1w=(wx_1x_2)^*\in L^1(\tau).\]
Furthermore, $w,x_1x_2\geq 0$ and so $w^{1/2}x_1x_2w^{1/2}\in L^{1}(\tau)$, by Proposition \ref{P3.4 DDP2}. Next we show that $x\leq 2k x_1 x_2$.  It follows from the first inequality in Proposition \ref{L4.8.16 BS}(2) and the assumption that $\varphi_2^{-1}\leq k \varphi_1^{-1}.\varphi_3^{-1}$ that 
\[t\leq \varphi_2^{-1}(t).(\varphi_2^*)^{-1}(t)\leq  k \varphi_1^{-1}(t).\varphi_3^{-1}(t). (\varphi_2^*)^{-1}(t).\]
Multiplying through by $(\varphi_3^*)^{-1}(t)$ and applying the second inequality in  Proposition \ref{L4.8.16 BS}(2) we therefore obtain 
\[t.(\varphi_3^*)^{-1}(t)\leq   k \varphi_1^{-1}(t).\varphi_3^{-1}(t).(\varphi_3^*)^{-1}(t). (\varphi_2^*)^{-1}(t) \leq 2kt \varphi_1^{-1}(t). (\varphi_2^*)^{-1}(t).\] 
Using $t= \varphi_3^*(s)$, we obtain 
\[ s \leq (\varphi_3^*)^{-1}\circ \varphi_3^*(s) \leq 2k \varphi_1^{-1}\circ \varphi_3^*(s).(\varphi_2^*)^{-1}\circ \varphi_3^*(s),\] for all $s>0$ and therefore $x\leq 2k x_1 x_2$, using the properties of the Borel functional calculus and the definitions of $x_1$ and $x_2$. It follows that $w^{1/2}xw^{1/2}\leq w^{1/2}2kx_1x_2w^{1/2}$, by \cite[Proposition 4.5]{key-Pag}, and therefore
\begin{eqnarray*}
\tau(w^{1/2}xw^{1/2})&\leq &\tau(w^{1/2}2kx_1x_2w^{1/2}) \qquad \text{using the positivity of $\tau$} \\
&=&2k\tau(wx_1x_2) \qquad \text{using Proposition \ref{P3.4 DDP2}} \\
&\leq& \tfrac{4k}{1-\epsilon}\norm{M_w}\norm{x}_{\varphi_3^*} \qquad \text{using (\ref{e2 T1 25/08/17}).}
\end{eqnarray*}
It follows that $(w^{1/2}x)w^{1/2}\in L^{1}(\tau)$. Since $x$, $x_1$ and $x_2$ commute, we can use the fact that $x\leq 2kx_1x_2$ to show that $|xw|^2\leq 4k^2|x_2x_1w|^2$ and hence $|xw|\leq 2k|x_2x_1w|$, since taking square roots is an operator-monotone function (see \cite[Corollary 3.2]{key-Hoa14}). It follows that $wx\in L^1(\tau)$. Since, apart from the restriction on the size of the norm, $x$ was an arbitrary element of $L^{\varphi_3^*}(\tau)^+$, It follows from K\"othe duality that $w\in L^{\varphi_3}(\tau)$. We can therefore apply Proposition \ref{P3.4 DDP2} to obtain $\tau(w^{1/2}xw^{1/2})=\tau(wx)$ and therefore, since the above holds for all $\epsilon >0$, we have
\begin{eqnarray}
|\Gamma(x)|=|\tau(wx)|=\tau(w^{1/2}xw^{1/2})\leq 4k\norm{M_w}\norm{x}_{\varphi_3^*}. \label{e10 T1 25/08/17}
\end{eqnarray}
Given $y\in L^{\varphi_3^*}(\tau)$, we may now use the polar decomposition $y=u|y|$ in terms of some partial isometry $u\in \mathcal{A}$, to conclude from the above that
\begin{eqnarray*}
|\Gamma(y)| &=& |\tau(wv|y|)| \\
& =&\tau((w^{1/2}v|y|^{1/2})(|y|^{1/2}w^{1/2}))\\
&\leq&\tau(w^{1/2}v|y|vw^{1/2})^{1/2}\tau(w^{1/2}|y|w^{1/2})^{1/2}\\
&\leq& (4k\norm{M_w}\norm{v|y|v}_{\varphi_3^*})^{1/2}(4k\norm{M_w}\norm{|y|}_{\varphi_3^*})^{1/2}\\
&\leq& 4k\norm{M_w}\norm{y}_{\varphi_3^*},
\end{eqnarray*}
where we used the fact that $\norm{v|y|v}_{\varphi_3^*}\leq \norm{|y|}_{\varphi_3^*}=\norm{y}_{\varphi_3^*}$ to obtain the final inequality. Clearly, $\|\Gamma\|\leq 4k\norm{M_w}$. Since Orlicz spaces are strongly symmetric spaces, we may apply \cite[Theorem 5.11]{key-DDP2} to conclude that $\norm{\Gamma}=\norm{w}_{\varphi_3}$ and hence that $\norm{w}_{\varphi_3}\leq 4k\norm{M_w}$.
\end{proof}
\end{theorem}

Applying this result to $L^p$-spaces, we obtain the following.

\begin{corollary}\label{C2 07/09/17}
Suppose $1< q < p$ and $r$ is such that $1/p+1/r=1/q$. If $w\in S(\mathcal{A},\tau)$, then $M_w$ is a bounded multiplication operator from $L^p(\tau)$ into $L^q(\tau)$ if and only if $w\in L^r(\tau)$.
\begin{proof}
If $\varphi_1(t):=t^p$, $\varphi_2(t)=t^q$ and $\varphi_3(t)=t^r$, then
\[\varphi^{-1}_1(t)\varphi_3^{-1}(t)=t^{1/p}t^{1/r}=t^{1/q} = \varphi_2^{-1}(t), \qquad \forall  t\geq 0.\] 
\end{proof}
\end{corollary}

In \cite{key-Lab14} sufficient conditions are obtained for the existence of multipliers between Orlicz spaces when the Orlicz functions are related by certain composition relations. We show that multiplication operators can be completely characterized and norm estimates obtained under similar circumstances.

\begin{theorem}\label{T1 25/08/17}
Suppose $w\in S(\mathcal{A},\tau)$ and $\psi, \varphi_1,\varphi_2$ are Orlicz functions. If $\varphi_3:=\psi^* \circ \varphi_2$ is an Orlicz function, $\psi \circ \varphi_2=\varphi_1$, and $\varphi_2$ is an Orlicz function satisfying the $\nabla'$-condition, then
\begin{enumerate}
\item there exists $0<c\leq 1$ so that $\varphi_2^{-1}(t) \leq c \varphi_1^{-1}(t) \varphi_3^{-1}(t)$ for all $t> 0$; 
\item $M_w$ is a bounded multiplication operator from $L^{\varphi_1}(\tau)$ into $L^{\varphi_2}(\tau)$ if and only if $w\in L^{\varphi_3}(\tau)$. Moreover, if this is the case, then \[\norm{M_w}\leq \frac{2}{c}\norm{w}_{\varphi_3}\leq 8\norm{M_w},\] 
\end{enumerate}
where $c>0$ is such that $\varphi_2^{-1}(st)\leq c \varphi_2^{-1}(s)\varphi_2^{-1}(t)$ for all $s,t \geq 0$ (see  Lemma \ref{L2 26/09/2017}). 
\begin{proof}
To prove (1), we start by showing that $\varphi_2^{-1}\circ \psi^{-1}\leq \varphi_1^{-1}$. Using the fact that $\varphi_1=\psi\circ \varphi_2$, we have that 
\[\varphi_1(\varphi_2^{-1}(t))=\psi\circ \varphi_2(\varphi_2^{-1}(t))\leq \psi(t),\] by Proposition \ref{Rp.276 BS} and the fact that $\psi$ is increasing. Replacing $t$  with $\psi^{-1}(t)$ in the inequality above we obtain 
\[\varphi_1\circ \varphi_2^{-1}(\psi^{-1}(t))\leq \psi(\psi^{-1}(t)) \leq t.\] 
Apply $\varphi_1^{-1}$ to the inequality above, and use the the fact that $\varphi_1^{-1}$ is increasing, to conclude from Proposition \ref{Rp.276 BS} that 
\[\varphi_2^{-1}(\psi^{-1}(t))\leq \varphi_1^{-1}\left(\varphi_1\circ \varphi_2^{-1}(\psi^{-1}(t))\right)\leq  \varphi_1^{-1}(t),\] 
as desired. A similar proof shows that $\varphi_2^{-1}\circ (\psi^*)^{-1}\leq \varphi_3^{-1}$.
 
Next, suppose $t>0$. Then using Proposition \ref{L4.8.16 BS} and the fact that $\varphi_2$ is increasing, we have that 
$\varphi_2^{-1}(t)\leq \varphi_2^{-1}\left(\psi^{-1}(t).(\psi^*)^{-1}(t)\right)$. On applying Lemma \ref{L2 26/09/2017}, it then follows 
that $\varphi_2^{-1}(t)\leq c \varphi_2^{-1}\left(\psi^{-1}(t)\right).\varphi_2^{-1}\left((\psi^*)^{-1}(t)\right)$. The inequalities 
verified in the first part of the proof now enable us to conclude that $\varphi_2^{-1}(t)\leq c\varphi_1^{-1}(t).\varphi_3^{-1}(t)$.

To prove (2), we start by noting that if $M_w$ is a bounded multiplication operator from $L^{\varphi_1}(\tau)$ into $L^{\varphi_2}(\tau)$, then it follows from the first part of this Theorem and Theorem \ref{T1 07/09/17}(2) that  $w\in L^{\varphi_3}(\tau)$ and $\norm{w}_{\varphi_3}\leq 4c \norm{M_w}$. 

It remains to prove the sufficiency and reverse inequality in (2). Given any $\lambda\geq 0$ and $0<\epsilon<1$, we may use convexity to conclude that $\epsilon^{-1}\varphi_2(\lambda)\leq  \varphi_2(\epsilon^{-1}\lambda)$, and hence that $\psi^*(\epsilon^{-1}\varphi_2(\lambda))\leq  \psi^*\circ\varphi_2(\epsilon^{-1}\lambda)=\varphi_3(\epsilon^{-1}\lambda)$. Given any positive Borel function $f$ for which $\varphi_3(\epsilon^{-1}f)$ is finite almost everywhere, it is clear that $\psi^*(\epsilon^{-1}\varphi_2(f))$ is then also finite almost everywhere with $\psi^*(\epsilon^{-1}\varphi_2(f))\leq\varphi_3(\epsilon^{-1}f)$. Since $\psi^*$ is an Orlicz function, this can clearly only be the case if $\varphi_2(f)$ itself is also finite almost everywhere. Given any $y\in L^{\varphi_3}(\tau)^+$ with $\norm{y}_{\varphi_3}<1$, it follows from Proposition \ref{P2.2 Lab13} that for any $\epsilon$ with $\|y\|_{\varphi_3}<\epsilon< 1$, we will have that $\varphi_3(\epsilon^{-1}y)\in S(\mathcal{A},\tau) $. On using the Borel functional calculus with $y$ playing the role of a Borel function, we may now use the above calculations to conclude that each of $\psi^*(\epsilon^{-1}\varphi_2(y))$ and $\varphi_2(y)$ are operators affiliated to $\mathcal{A}$ (see \cite[Lemma 9.4.7 and Theorem 9.4.8]{key-Tak79}), for which we have that $\psi^*(\epsilon^{-1}\varphi_2(f))\leq\varphi_3(\epsilon^{-1}f)$.  We may then use Proposition 
\ref{L4.8.16 BS} to conclude that $\epsilon^{-1}\varphi_2(y)\leq (\psi^*)^{-1}(\psi^*(\epsilon^{-1}\varphi_2(y))) \leq (\psi^*)^{-1}(\varphi_3(\epsilon^{-1}y))$. Since by Lemma \ref{L1 20/09/17}, $(\psi^*)^{-1}(\varphi_3(\epsilon^{-1}y)) \in S(\mathcal{A},\tau)$, the 
preceding inequality ensures that also $\varphi_2(y)\in S(\mathcal{A},\tau)$. Given that the specific choice of $\epsilon$ and the preceding 
inequalities ensure that $\tau(\psi^*(\epsilon^{-1}\varphi_2(y)))\leq\tau(\varphi_3(\epsilon^{-1}y))\leq 1$, it follows that 
$\varphi_2(y)\in L^{\psi^*}(\tau)^+$ with $\norm{\varphi_2(y)}_{\psi^*}\leq \norm{y}_{\varphi_3}$. One may similarly show that if 
$v\in L^{\varphi_1}(\tau)^+$ with $\|v\|_{\varphi_1}<1$, then $\varphi_2(v)\in L^{\psi}(\tau)^+$ with $\norm{\varphi_2(v)}_{\psi}\leq \norm{v}_{\varphi_1}$. 

Recall that by assumption there exists $0< c\leq 1$ such that $\varphi_2(cst)\leq \varphi_2(s)\varphi_2(t)$ for all $s,t\geq 0$.
Let $n\in \mathbb{N}$ be given with $n>1$, and suppose we are given $w\in L^{\varphi_3}(\tau)$ and $x\in L^{\varphi_1}(\tau)$ with $\norm{w}_{\varphi_3}=\tfrac{cn}{n+1}$ and $\norm{x}_{\varphi_1}=\tfrac{n+1}{2n}$. Both have norm less than 1, so $\varphi_2(w)\in L^{\psi^*}(\tau)^+$ with $\norm{\varphi_2(w)}_{\psi^*}\leq \norm{w}_{\varphi_3}$, and $\varphi_2(x)\in L^{\psi}(\tau)^+$ with $\norm{\varphi_2(x)}_{\psi}\leq \norm{x}_{\varphi_1}$. Then 
\begin{eqnarray}
 \tau(\varphi_2(|wx|))=\int_0^\infty \varphi_2(\svft{wx}{t})dt, \label{e1 T1 25/08/17}
 \end{eqnarray}
 by Lemma \ref{L2.1 Lab13}. Furthermore, $\varphi_2$ is increasing and it is easily checked that $t \mapsto \varphi_2(e^t)$ is convex. Therefore
 \begin{eqnarray} 
&&\int_0^\infty \varphi_2(\svft{wx}{t})dt \leq \int_0^\infty \varphi_2(\svft{w}{t}\svft{x}{t})dt \leq \int_0^\infty \varphi_2(\svft{w/c}{t})\varphi_2(\svft{x}{t})dt, \label{e1b T1 25/08/17}
\end{eqnarray}
by Proposition \ref{T4.2 Fack86}(3). It follows by what has been shown already that $\varphi_2(\svf{w/c})\in L^{\psi^*}(0,\infty)$ and $\varphi_2(\svf{x})\in L^{\psi}(0,\infty)$. Therefore,  by Proposition \ref{P2.3 Lab13}, 
\begin{eqnarray*}
\int_0^\infty \varphi_2(\svft{w/c}{t})\varphi_2(\svft{x}{t})dt&\leq& \norm{\varphi_2(\svf{w/c})}_{L^{\psi^*}(0,\infty)}^0 \norm{\varphi_2(\svf{x})}_{L^\psi(0,\infty)} \\
&\leq&2\norm{\varphi_2(\svf{w/c})}_{L^{\psi^*}(0,\infty)} \norm{\varphi_2(\svf{x})}_{L^\psi(0,\infty)}, 
\end{eqnarray*} 
where the second inequality follows by \cite[Theorem 4.8.14]{key-BS88}. Combining this with (\ref{e1 T1 25/08/17}) and (\ref{e1b T1 25/08/17}) we therefore obtain
\begin{eqnarray*}
\tau(\varphi_2(|wx|))&\leq& 2\norm{\varphi_2(\svf{w/c})}_{L^{\psi^*}(0,\infty)} \norm{\varphi_2(\svf{x})}_{L^\psi(0,\infty)}\\
&=& 2\norm{\svf{w/c}}_{L^{\varphi_3}(0,\infty)} \norm{\svf{x}}_{L^{\varphi_1}(0,\infty)} =1, 
\end{eqnarray*}
since $\norm{\svf{w/c}}_{L^{\varphi_3}(0,\infty)}=\norm{w/c}_{\varphi_3}=\frac{n}{n+1}$ and similarly $\norm{\svf{x}}_{L^{\varphi_1}(0,\infty)}=\tfrac{n+1}{2n}$. It follows that \[\norm{wx}_{\varphi_2}=\inf\{\lambda>0:\tau(\varphi_2(|wx|/\lambda))\leq 1\}\leq 1 = 2 \norm{w/c}_{\varphi_3}\norm{x}_{\varphi_1}.\] 
It follows that $M_w$ is a bounded multiplication operator from $L^{\varphi_1}(\tau)$ into $L^{\varphi_2}(\tau)$ and $\norm{M_w}\leq \frac{2}{c}\norm{w}_{\varphi_3}$. 
\end{proof}
\end{theorem}

We finish this subsection by showing that the previous result also applies to $L^p$-spaces.

\begin{corollary}\label{CT1 25/08/17}
Suppose $1<q<p$ and let $r>1$ be such that $1/p+1/r=1/q$.  If $w\in S(\mathcal{A},\tau)$, then $M_w$ is a bounded multiplication operator from $L^p(\tau)$ into $L^q(\tau)$ if and only if $w\in L^r(\tau)$.
\begin{proof}
 Let $\varphi_1(t)=t^p$, $\varphi_2(t)=t^q$ and $\psi(t)=t^{p/q}$. Then $\varphi_1$, $\varphi_3$ and $\psi$ are Orlicz functions. Furthermore, $\psi \circ \varphi_2 (t) = (t^q)^{p/q}=t^p=\varphi_1(t)$ and $r=qp/(p-q)$. If we let $(p/q)'$ denote the conjugate exponent of $p/q$, then it is easily checked that $(p/q)'=p/(p-q)=r/q$. A straightforward calculation shows that $\psi^*(t)=\frac{1}{(p/q)^{(r/q)/p}(r/q)}t^{r/q}$. If we let $\varphi_3=\psi^*\circ \varphi_2$, then $\varphi_3=ct^r$, where $c=\frac{1}{(p/q)^{(r/q)/p}(r/q)}t^{r/q}$. It follows that $\varphi_3$ is an Orlicz function and $L^{\varphi_3}(\tau)=L^r(\tau)$, with $\norm{x}_{\varphi_3}=c^{1/p}\norm{x}_r$ for every $x\in L^r(\tau)$. Furthermore, $\varphi_2$ is an Orlicz function satisfying the $\nabla'$-condition in that 
 \[\varphi_2^{-1}(st)=(st)^{1/q}=s^{1/q}t^{1/q}=\varphi_2^{-1}(s)\varphi_2^{-1}(t) \qquad \forall s,t\geq 0.\]
The result therefore follows by Theorem \ref{T1 25/08/17}.
\end{proof}
\end{corollary}

\subsection{Multiplication operators on $L^p$-spaces}

In this subsection we consider multiplication operators from $L^p(\tau)$ into $L^q(\tau)$. It follows from \cite[Proposition 5]{key-Han15} that $M_w$ is a bounded multiplication operator from $L^p(\tau)$ ($1\leq p \leq \infty$) into itself if and only if $w\in \mathcal{A}$, in which case $\norm{M_w}=\norminf{w}$.  In Corollary \ref{C2 07/09/17} and Corollary \ref{CT1 25/08/17} we used the theory for Orlicz spaces to conclude that if $1 < q<p<\infty$ and $w\in S(\mathcal{A},\tau)^+$, then $M_w$ is a bounded multiplication operator from $L^p(\tau)$ into $L^q(\tau)$ if and only if $w\in L^{r}(\tau)$, where $1/p+1/r=1/q$ (this result is also claimed in \cite[Example 1(ii)]{key-Han15} although no proof is given). In this subsection we will see that a direct proof will however enable us to determine the norm of the multiplication operator exactly in this case. We also consider the case $1\leq p<q<\infty$.

\begin{theorem}\label{T1 30/05/17}
Suppose $1< q<p<\infty$ and $w\in S(\mathcal{A},\tau)^+$. Then $M_w$ is a bounded multiplication operator from $L^p(\tau)$ into $L^q(\tau)$ if and only if $w\in L^{r}(\tau)$, where $1/p+1/r=1/q$. Furthermore, if this is the case, then $\norm{M_w}=\normx{w}{r}$.
\begin{proof}
Suppose $w\in L^r(\tau)$. Then, using Proposition \ref{T4.2 Fack86}(2), we obtain \[\normx{M_wx}{q}=\normx{wx}{q}\leq \normx{w}{r}\normx{x}{p}.\]
It follows that $M_w$ is a bounded multiplication operator from $L^p(\tau)$ into $L^q(\tau)$ and $\norm{M_w}\leq \normx{w}{r}$. 

Conversely, suppose $M_w$ is a bounded multiplication operator from $L^p(\tau)$ into $L^q(\tau)$. Let $x\in L^{r'}(\tau)$ be given, where $r'$ is the conjugate index to $r$. If $x=u|x|$ is the polar form of $x$, we set $x_p=u|x|^{r'/p}$ and $x_{q'}=|x|^{r'/q'}$. It is an easy exercise to see that $x_p\in L^p(\tau)$ and $x_{q'}\in L^{q'}(\tau)$. By hypothesis, $wx_p\in L^q(\tau)$, and so $wx=wx_px_{q'}\in L^1(\tau)$. K\"othe duality now ensures that $w\in L^r(\tau)$. 

We proceed to prove the equality of $\norm{M_w}$ and $\normx{w}{r}$. It is not difficult to conclude from the fact that $w\in L^r(\tau)^+$, that $w^{r/p}\in L^p(\tau)$ with $\|w^{r/p}\|_{p}=(\|w\|_{r})^{r/p}$. But then
$$\|M_ww^{r/p}\|_{q}=\|w^{r/q}\|_{q}=(\|w\|_{r})^{r/q}=\|w\|_{r}.(\|w\|_{r})^{r/p}=\|w\|_{r}.\|w^{r/p}\|_{p}.$$This clearly ensures that $\|M_w\|\geq \|w\|_{r}$, and hence that equality of norms must hold.
\end{proof}
\end{theorem}

Next, we consider the case $p<q$.

\begin{remark}\label{R2 12/07/18}
It is claimed in \cite[Example 1(i)]{key-Han15} that the space of multipliers from $L^p(\tau)$ into $L^q(\tau)$ ($p<q$) consists just of the the zero operator. This cannot be true in general as Example \ref{E1 13/07/18} shows. 
\end{remark}

 We start by showing that in this setting it suffices to consider purely atomic von Neumann algebras. 

\begin{theorem}\label{T2 23/10/17}
Suppose $1\leq p<q<\infty$ and $c$ is the central projection such that $c\mathcal{A}$ is atomic and $c^{\perp}\mathcal{A}$ is non-atomic. If $w\in S(\mathcal{A},\tau)^+$ is such that $M_w$ is a bounded multiplication operator from $L^p(\tau)$ into $L^q(\tau)$, then $wc^{\perp}=0$.
\begin{proof}
Suppose that $M_w$ is a bounded multiplication operator from $L^p(\tau)$ into $L^q(\tau)$ and assume that $e=e^{wc^\perp}(\lambda,\infty)$ is non-zero for some $\lambda>0$. Since $e\in c^\perp\mathcal{A}$ and $c^\perp\mathcal{A}$ is non-atomic, it follows that given $0<\alpha<\tau(e)$, we may select a sequence $(e_n)_{n=1}^\infty\subset c^\perp\mathcal{A}$ of mutually orthogonal subprojections of $e$ such that $\tau(e_n)=\alpha/2^n$. Let $v_n:= e_n/\tau(e_n)^{1/q}$. Since $1/p>1/q$ and for each $n$ we have that
$$\tau(|v_n|^p)^{1/p}=\tau(e_n)^{1/p-1/q}=\left(\alpha/2^n\right)^{1/p-1/q},$$it is clear that $(v_n)_{n=1}^\infty$ is a sequence which converges to zero in $L^p(\tau)$. So by continuity of $M_w$, $(wc^\perp v_n)_{n=1}^\infty=(w v_n)_{n=1}^\infty$ must converge to zero in $L^q(\tau)$. But this cannot be, since the fact that $e_n\leq e$, ensures that for all $n$ we have that $wc^\perp v_n= \frac{wc^\perp e_n}{\tau(e_n)^{1/q}}\geq \lambda\frac{e_n}{\tau(e_n)^{1/q}}$, with $\|wc^\perp v_n\|_{q}\geq \lambda\tau(\frac{e_n}{\tau(e_n)})^{1/q}=\lambda$. This clear contradiction establishes the claim.
\end{proof}
\end{theorem}

Since we are dealing with atomic von Neumann algebras, it suffices to consider a (possibly uncountable) direct sum of (possibly infinite) factors. We will therefore consider the situation on each such factor before investigating the general case. We start by proving a simple lemma that will help us in this regard. 

\begin{lemma}\label{L2 20/10/17}
Suppose $\mathcal{A}=\mathcal{B}(H)$, $\tau$ is a faithful, semi-finite normal trace on $\mathcal{A}$ and $1\leq p<\infty$. If $w\in L^p(\tau)$ and $e$ is a projection onto a one-dimensional subspace of $H$, then \[\norm{we}_{L^p(\tau)}=\norminf{we}\tau(e)^{1/p}.\]
\begin{proof}
 Note that
 \[|we|^2=e|w|^2 e = \lambda e,\]
  for some $\lambda\geq 0$, by \cite[Proposition 6.4.3]{key-K2}. It follows that for any $\xi\in H$, 
  \[\norm{(we)\xi}^2=\ip{we \xi}{we \xi}= \ip{\lambda e \xi}{\xi}=\norm{(\lambda^{1/2}e)\xi}^2\] and so $\norminf{we}=\lambda^{1/2}$.  Therefore 
 \begin{eqnarray*}
 \tau(|we|^p)^{1/p}= \tau((|we|^2)^{p/2})^{1/p} =\tau(\lambda^{p/2} e)^{1/p} =\lambda^{1/2} \tau(e)^{1/p} =\norminf{we}\tau(e)^{1/p}. 
 \end{eqnarray*}
\end{proof}
\end{lemma}

For multiplication operators between $L^p$-spaces associated with factors we have the following characterization. 

\begin{theorem}\label{T1b 13/10/17}
Suppose $\mathcal{A}=\mathcal{B}(H)$, $\tau$ is a faithful semi-finite normal trace on $\mathcal{A}$ and $1\leq p<q<\infty$. Then for all $w\in S(\mathcal{A},\tau)$, $M_w$ is a bounded multiplication operator from $L^p(\tau)$ into $L^q(\tau)$ (note that in this case $S(\mathcal{A},\tau)=\mathcal{A}$). Furthermore,  
\[\norm{M_w}= k^{-1/s}\norminf{w},\]
where $k$ is the trace of any projection onto a one-dimensional subspace of $H$ and $s>0$ is such that $1/q+1/s=1/p$.
\begin{proof}
It follows from \cite[Propositions 8.5.3 \& 8.5.5]{key-K2} that $\tau(\cdot)=k\,tr(\cdot)$, for some $k>0$, where $tr(\cdot)$ denotes the canonical trace on $\mathcal{B}(H)$. Since all projections onto one-dimensional subspaces of $H$ have the same trace, $k=\tau(e)$, where $e$ is any such projection. We already noted that in this case $S(\mathcal{A},\tau)=\mathcal{A}$, since $\mathcal{A}=\mathcal{B}(H)$. Furthermore, if $1\leq p<q<\infty$, then $L^p(tr)$ is continuously embedded into $L^q(tr)$ and
\[ tr(|x|^q)^{1/q} \leq tr(|x|^p)^{1/p}\] for any $x\in L^p(tr)$, by \cite[Proposition 4.5]{key-Diestel95}. It follows that
\begin{eqnarray}
\norm{x}_{L^q(\tau)} \leq k^{-1/s}\norm{x}_{L^p(\tau)} \label{eP1b 20/10/17} \qquad \forall x \in L^p(\tau),
\end{eqnarray} 
where $s>0$ is such that $1/q+1/s=1/p$. If we let $w\in S(\mathcal{A},\tau)= \mathcal{A}$, then for any $x\in L^p(\tau)$ we have  
\begin{eqnarray*}
 \norm{wx}_{L^q(\tau)}\leq \norminf{w}\norm{x}_{L^q(\tau)} \leq \norminf{w}k^{-1/s}\norm{x}_{L^p(\tau)},  
\end{eqnarray*}
using (\ref{eP1b 20/10/17}). It follows that $M_w$ is a bounded multiplication operator from $L^p(\tau)$ into $L^q(\tau)$ and 
\[\norm{M_w}\leq k^{-1/s} \norminf{w}.\]

To prove the reverse inequality, let $\{p_\alpha\}_{\alpha\in \mathbb{A}}$ denote the collection of all distinct projections onto one dimensional subspaces of $H$ and consider $x_\alpha:=\frac{p_\alpha}{\tau(p_\alpha)^{1/p}}$. Then $x_\alpha\in L^p(\tau)$ and $\norm{x_\alpha}_{L^p(\tau)}=1$, for each 
$\alpha\in \mathbb{A}$. Therefore 
 \begin{eqnarray}
 &&\norm{M_w}\geq  \norm{M_wx_\alpha}_{L^q(\tau)} =\tau(p_\alpha)^{-1/p}\norminf{wp_\alpha}\tau(p_\alpha)^{1/q} =\frac{\norminf{wp_\alpha}}{\tau(p_\alpha)^{1/s}}, \label{e4b T1 13/10/17}
 \end{eqnarray}
 by Lemma \ref{L2 20/10/17}. Since (\ref{e4b T1 13/10/17})  holds for all $\alpha \in \mathbb{A}$ and $\tau(p_\alpha)=k$, for each $\alpha$, we have that 
\begin{eqnarray}
k^{-1/s}\supu{\alpha \in \mathbb{A}}\norminf{wp_\alpha}\leq \norm{M_w}<\infty. \label{e2b T1 23/10/17}
\end{eqnarray} 
Since $\norminf{w}=\sup\{\norm{w\xi}:\norm{\xi}=1\}$, there exists a sequence $\{\xi_n\}_{n=1}^\infty$ of norm one elements such that $\norm{w \xi_n}> \norminf{w}-1/n$, for every $n\in \mathbb{N}$. Let $p_{n}$ denote the projection onto the one-dimensional subspace generated by $\xi_n$. Then 
\[\norminf{wp_{n}} \geq \norm{wp_{n}\xi_n}=\norm{w \xi_n}>\norminf{w}-1/n.\] 
It follows that $\supu{\alpha \in \mathbb{A}}\norminf{wp_{\alpha}}=\norminf{w}$ and therefore by inequality (\ref{e2b T1 23/10/17}) that 
\[k^{-1/s}\norminf{w}\leq \norm{M_w}.\]
\end{proof}
\end{theorem}

We will need the following lemma in order to move from a factor to a direct sum of factors.

\begin{lemma}\label{LP1 10/10/17}
If $\{p_\alpha\}_{\alpha \in \mathbb{A}}$ is a (possibly uncountable) family of mutually orthogonal central projections, $q\geq 1$ and $x\in L^q(\tau)$, then $|x\summ{\alpha \in \mathbb{A}}{}p_\alpha|^q=\summ{\alpha\in \mathbb{A}}{}|xp_\alpha|^q$ and $\normx{x\summ{\alpha \in \mathbb{A}}{}p_\alpha}{q}^q=\summ{\alpha \in \mathbb{A}}{}\tau(|xp_\alpha|^q)$.
\begin{proof}
Suppose $x\in L^q(\tau)$ and $c$ is a central projection. For any Borel function we know from \cite[Lemma 5.6.31]{key-K1} that $g(|x|e)e=g(|x|)e$. But then also $g(|x|e)e^\perp=g((|x|e)e^\perp)e^\perp=g(0)e^\perp$. For the specific function $g(t)=t^q$, these facts ensure that $|x|^q p_{\alpha}=|xp_\alpha|^q$ for each $\alpha$. Suppose $\mathbb{B}$ is a finite subcollection of $\mathbb{A}$. Since $\{p_\alpha\}_{\alpha \in \mathbb{B}}$ is a collection of mutually orthogonal central projections (and hence $\summ{\alpha \in \mathbb{B}}{}p_\alpha$ is also a central projection), this enables us to conclude that
\begin{equation}
|x\summ{\alpha \in \mathbb{B}}{}p_\alpha|^q=|x|^q\summ{\alpha \in \mathbb{B}}{}p_\alpha=\summ{\alpha\in \mathbb{B}}{}|x|^qp_\alpha=\summ{\alpha \in \mathbb{B}}{}|xp_\alpha|^q. \label{e1c 10/10/17}
\end{equation}
Furthermore, 
\begin{eqnarray}
&&|x\summ{\alpha \in \mathbb{A}}{}p_\alpha|=|x|^{1/2}\left(\summ{\alpha \in \mathbb{A}}{}p_\alpha\right)|x|^{1/2}\geq |x|^{1/2}\left(\summ{\alpha \in \mathbb{B}}{}p_\alpha\right)|x|^{1/2}=|x\summ{\alpha \in \mathbb{B}}{}p_\alpha|.  \label{e1d 10/10/17}
\end{eqnarray} 
Using (\ref{e1c 10/10/17}) and (\ref{e1d 10/10/17}), it follows that
\[\summ{\alpha \in \mathbb{B}}{}\tau(|xp_\alpha|^q)=\normx{x\summ{\alpha \in \mathbb{B}}{}p_\alpha}{q}^q \leq \normx{x\summ{\alpha \in \mathbb{A}}{}p_\alpha}{q}^q. \]
Since this holds for every finite subcollection of $\mathbb{A}$, this ensures that $\summ{\alpha \in \mathbb{A}}{}\tau(|xp_\alpha|^q)$ converges, with $\summ{\alpha \in \mathbb{A}}{}\tau(|xp_\alpha|^q)\leq \normx{x\summ{\alpha \in \mathbb{A}}{}p_\alpha}{q}^q < \infty$. Therefore $\tau(|xp_\alpha|^p)$, and hence also $xp_\alpha$, is non-zero for at most countably many $\alpha \in \mathbb{A}$. Let $(p_n)_{n=1}^{\infty}$ denote the collection of projections for which this holds. 

Since $\summ{n=1}{K}p_n\uparrow \summ{n=1}{\infty}p_n$ and these are all central projections, we may use (\ref{e1c 10/10/17}) to conclude that  $\summ{n=1}{K}|xp_n|^q\uparrow  \left|x\left(\summ{n=1}{\infty}p_n\right)\right|^{q}$ and therefore 
\[\summ{\alpha\in \mathbb{A}}{}|xp_\alpha|^q=\summ{n=1}{\infty}|xp_n|^q=\left|x\left(\summ{n=1}{\infty}p_n\right)\right|^{q}=\left|x\summ{\alpha \in \mathbb{A}}{}p_\alpha\right|^q,\] 
as desired. Furthermore, 
\begin{eqnarray}
&&\summ{n=1}{K}\tau(|xp_n|^q)=\tau(\summ{n=1}{K}|xp_n|^q)\uparrow \tau( |\left(\summ{n=1}{\infty}p_n\right)x|^{q}) =\normx{\left(\summ{n=1}{\infty}p_n\right)x}{q}^q. \label{e3 10/10/17}
\end{eqnarray} 
Since these are the only projections for which $xp_\alpha\neq 0$, it follows that $\normx{x\summ{\alpha \in \mathbb{A}}{}p_\alpha}{q}^q=\summ{\alpha \in \mathbb{A}}{}\tau(|xp_\alpha|^q)$.
\end{proof}
\end{lemma}

We are now in a position to characterize multiplication operators from an $L^p$-space into an $L^q$-space for the case $p<q$. 

\begin{theorem}\label{T3 26/10/17}
Suppose $(\mathcal{A},\tau)$ is a semi-finite von Neumann algebra, $1\leq p<q<\infty$, $w\in S(\mathcal{A},\tau)^+$ and $c$ is the central projection such that $c \mathcal{A}$ is atomic (i.e. $c\mathcal{A} \cong \dsumx{\alpha \in \mathbb{A}}{} \mathcal{B}(H_\alpha)$) and $c^{\perp}\mathcal{A}$ is non-atomic. Let $p_\alpha$ denote the central projection such that $p_\alpha \mathcal{A} \cong \mathcal{B}(H_\alpha)$ and let $k_\alpha$ denote the trace of a projection in $\mathcal{B}(H_\alpha)$ onto a one-dimensional subspace of $H_\alpha$. Then $M_w$ is a bounded multiplication operator from $L^p(\tau)$ into $L^q(\tau)$ if and only if $wc=w$ and $\supu{\alpha \in \mathbb{A}} \frac{\norminf{wp_\alpha}}{k_\alpha^{1/s}}<\infty$, where $s>0$ is such that $1/q+1/s=1/p$. Furthermore, if this is the case, then 
\[\norm{M_w}= \supu{\alpha \in \mathbb{A}} \frac{\norminf{wp_\alpha}}{k_\alpha^{1/s}}.\]
\begin{proof}
Suppose $wc=w$ and $\supu{\alpha \in \mathbb{A}} \frac{\norminf{wp_\alpha}}{k_\alpha^{1/s}}<\infty$. If $x\in L^p(\tau)$ with $\normx{x}{p}=1$, then  
\begin{eqnarray}
\normx{xp_\alpha}{p}^p=\left(\normx{xp_\alpha}{p}^{q}\right)^{p/q}\geq \normx{xp_\alpha}{p}^q, \label{e1 T1 23/10/17}
\end{eqnarray}
since $p/q<1$ and $\normx{xp_\alpha}{p}^q\leq \normx{x}{p}^q=1$. Furthermore, 
\begin{eqnarray*}
\normx{wx}{q}^q&=&\tau(|wx\summ{\alpha\in \mathbb{A}}{}p_\alpha|^q) \qquad \text{since $wc=w$}\\
&=&\summ{\alpha\in \mathbb{A}}{}\tau(|wxp_\alpha|^q) \qquad \text{by Lemma \ref{LP1 10/10/17}} \\
&=&\summ{\alpha\in \mathbb{A}}{}\normx{wxp_\alpha}{q}^q \\
&\leq&\summ{\alpha\in \mathbb{A}}{}k^{-q/s}\norminf{wp_\alpha}^q\normx{xp_\alpha}{p}^q \qquad \text{by Theorem \ref{T1b 13/10/17}} \\
&\leq&\supu{\alpha \in \mathbb{A}} \frac{\norminf{wp_\alpha}^q}{k_\alpha^{q/s}}\summ{\alpha\in \mathbb{A}}{}\normx{xp_\alpha}{p}^q \\
&\leq&\supu{\alpha \in \mathbb{A}} \frac{\norminf{wp_\alpha}^q}{k_\alpha^{q/s}}\summ{\alpha\in \mathbb{A}}{}\normx{xp_\alpha}{p}^p \qquad \text{using (\ref{e1 T1 23/10/17})} \\
&=&\supu{\alpha \in \mathbb{A}} \frac{\norminf{wp_\alpha}^q}{k_\alpha^{q/s}} \normx{x}{p}^p \qquad \text{by Lemma \ref{LP1 10/10/17}.} 
\end{eqnarray*}
It follows that 
\[\normx{wx}{q}\leq \supu{\alpha \in \mathbb{A}} \frac{\norminf{wp_\alpha}}{k_\alpha^{1/s}} \normx{x}{p}^{p/q}=\supu{\alpha \in \mathbb{A}} \frac{\norminf{wp_\alpha}}{k_\alpha^{1/s}}\normx{x}{p},\]
since $\normx{x}{p}=1$. Therefore, $M_w$ is a bounded multiplication operator from $L^p(\tau)$ into $L^q(\tau)$ and 
\[\norm{M_w}\leq \supu{\alpha \in \mathbb{A}} \frac{\norminf{wp_\alpha}}{k_\alpha^{1/s}}.\]
Conversely, suppose $M_w$ is a bounded multiplication operator from $L^p(\tau)$ into $L^q(\tau)$. Then $wc=w$, by Theorem \ref{T2 23/10/17}. Let $\tau_\alpha$ denote the restriction of $\tau$ to $p_\alpha \mathcal{A}$. Since $p_\alpha L^p(\tau)=L^p(p_\alpha\mathcal{A}, \tau_\alpha)$ and the action of $w$ on $L^p(p_\alpha\mathcal{A}, \tau_\alpha)$ is induced by $wp_\alpha$, we have that $M_{wp_\alpha}$ is a bounded multiplication operator from $L^p(p_\alpha \mathcal{A},\tau_\alpha)$ into $L^{q}(p_\alpha \mathcal{A},\tau_\alpha)$, for each $\alpha$ and 
\[\normx{wp_\alpha x}{q}\leq \norminf{p_\alpha}\normx{wx}{q}\leq \norm{M_w}\normx{x}{p},\] for each $x\in L^p(\tau)$. Using Theorem \ref{T1b 13/10/17}, it follows that $k_\alpha^{-1/s}\norminf{wp_\alpha} \leq \norm{M_{wp_\alpha}}\leq \norm{M_w}$, for each $\alpha$. Since this holds for all $\alpha$, we have that
\[\supu{\alpha \in \mathbb{A}} \frac{\norminf{wp_\alpha}}{k_\alpha^{1/s}}\leq \norm{M_w}.\]
\end{proof}
\end{theorem}

\section{Compactness of multiplication operators}

The characterizations and norm estimates obtained in the previous section will enable us to obtain characterizations of compactness in the same settings. It is easily checked that $w\in S(\mathcal{A},\tau)$ induces a compact multiplication operator between symmetric spaces $E$ and $F$ if and only if $|w|$ induces a compact multiplication operator. As in the previous section, it therefore suffices to consider positive elements in $S(\mathcal{A},\tau)$. We start by quoting a necessary condition for the compactness of multiplication operators which will be used throughout.

\begin{theorem} \cite[Theorem 4.2]{key-Lab14} \label{T4.2 Lab14}
Given two Orlicz functions $\psi_1$ and $\psi_2$ and $y\in S(\mathcal{A},\tau)$ such that $M_y:L^{\psi_1}(\tau) \rightarrow L^{\psi_2}(\tau)$ is compact. Then there exists a central projection $\tilde{c}$ such that $y\tilde{c}=y$ with $\tilde{c}\mathcal{A}$ being a direct sum of countably many finite type I factors.
\end{theorem}

The techniques employed to prove Theorem \ref{T4.2 Lab14} can easily be adapted to prove the same result for any pair of Banach function spaces which are intermediate spaces of the Banach couple $(L^\infty(\tau),L^1(\tau))$.

As with the boundedness of multiplication operators, a characterization for the compactness of endomorphic multiplication operators can be obtained in the general setting of symmetric spaces. We know from \cite[Proposition 5]{key-Han15} that if $E$ is a symmetric space and $w\in S(\mathcal{A},\tau)$, then $M_w$ is a bounded multiplication operator from $E$ into itself if and only if $w\in \mathcal{A}$. In considering the compactness of multiplication operators in the endomorphic setting it therefore suffices to consider $w\in \mathcal{A}$.

\begin{lemma} \label{L1 13/11/17} 
Suppose $E\subseteq S(\mathcal{A},\tau)$ is a symmetric space and $w\in \mathcal{A}^+$. Then $M_w$ is compact if and only if 
\[Z_\epsilon^w:=\{e^w(\epsilon,\infty)x:x\in E\}\]
is finite-dimensional for every $\epsilon>0$.
\begin{proof}
Suppose $M_w$ is compact. Let $x\in E$ and $\epsilon>0$. Then $wx\in E$, since $M_w$ is a bounded operator from $E$ into itself. Furthermore, $w$ and $e^{w}(\epsilon,\infty)$ commute and so 
\[M_w(e^w(\epsilon,\infty)x)=we^w(\epsilon,\infty)x=e^w(\epsilon,\infty)wx\in Z_\epsilon^w.\]
It follows that $Z_\epsilon^w$ is invariant under $M_w$. The restriction $\widetilde{M}_w$ of $M_w$ to $Z_\epsilon^w$ is therefore a compact operator from $Z_\epsilon^w$ into itself. Let $p=e^w(\epsilon,\infty)$. On passing to the reduced space $\mathcal{A}_p$ (see \cite[p.211]{key-Dodds14}), the inequality $pwp \geq \epsilon p$, ensures that as an element of $\mathcal{A}_p$, $pwp$ is strictly positive, and hence that there exists $v\in\mathcal{A}_p^+$ with $v\leq \tfrac{1}{\epsilon} p$, such that $pwpv=vpwp=p$. Since $p$ and $w$ commute, this will in $\mathcal{A}$ reduce to the statement that $wv=vw=p$. Using this element $v$, it is easily checked that $\widetilde{M}_w$ is invertible with the inverse given by the restriction of ${M}_v$ to $Z_\epsilon^w$. Since $\widetilde{M}_w$ is also compact, it therefore follows that  $Z_\epsilon^w$ is finite-dimensional.

Conversely, if $Z_\epsilon^w$ is finite dimensional for every $\epsilon>0$, then in particular $Z_{1/n}^w$ is finite dimensional for every $n\in \mathbb{N}^+$. Let $w_n$ be defined as $w_n=we^w(1/n,\infty)$. Then it is easily checked that $M_{w_n}(E)\subseteq Z_{1/n}^w$. It follows that  $M_{w_n}$ is finite rank for every $n\in \mathbb{N}^+$. Furthermore, if $x\in E$, then 
\begin{eqnarray*}
\norm{(M_w-M_{w_n})(x)}_E&=&\norm{we^w[0,1/n]x}_E \\
&\leq&  \norminf{we^w[0,1/n]}\norm{x}_E \\
&\leq& \tfrac{1}{n}\norm{x}_E.
\end{eqnarray*}
It follows that $\norm{M_w-M_{w_n}}\leq 1/n$ and hence that as the limit of a sequence of finite rank operators, $M_w$ is compact.
\end{proof}
\end{lemma}

\begin{theorem} \label{T2 13/11/17}
Let $E\subseteq S(\mathcal{A},\tau)$ be a symmetric space which is an intermediate space for the Banach couple $(L^\infty(\tau),L^1(\tau))$ and let $w\in \mathcal{A}^+$. Then $M_w$ is a compact multiplication operator from $E$ into itself if and only if 
\begin{itemize} 
\item there exists a sequence $(p_n)_{n=1}^\infty$ of mutually orthogonal central projections such that $w\tilde{c}=w$, where $\tilde{c}=\summ{n=1}{\infty} p_n$, $\tilde{c}\mathcal{A}=\dsumx{n=1}{\infty} \mathcal{A}_n$ and each $p_n\mathcal{A}=\mathcal{A}_n$ is a finite type 1 factor;
\item and $\norminf{wp_n}\rightarrow 0$ as $n \rightarrow 0$.
\end{itemize}
\begin{proof}
Suppose $M_w$ is compact. By Theorem \ref{T4.2 Lab14} and the comments following it, there exists a central projection $\tilde{c}$ with the desired properties. Since $p_n\mathcal{A}=\mathcal{A}_n$ is a finite type 1 factor, we can write $wp_n=\summ{k=N_n}{N_{n+1}-1}\lambda_{k}q_k$, where the $\lambda_{k}$'s ($N_n\leq k <N_{n+1}$) are the eigenvalues of $wp_n$, repeated according to multiplicity, and the $q_k$'s are mutually orthogonal projections onto one-dimensional subspaces of the eigenspace of $\lambda_{k}$. Recall that 
\[Z_\epsilon^w:=\{e^w(\epsilon,\infty)x:x\in E\}.\]
Fix $\epsilon>0$. If $\lambda_{k}>\epsilon$, then $q_k\leq e^w(\epsilon,\infty)$ (since $wq_k=\lambda_k q_k$) and so $q_k=e^w(\epsilon,\infty)q_k\in Z_\epsilon^w$. Since $Z_\epsilon^w$ is finite-dimensional, by Lemma \ref{L1 13/11/17}, there can only be a finite number of $\lambda_{k}$'s such that $\lambda_{k}>\epsilon$. It follows that $\lambda_k \rightarrow 0$ as $k \rightarrow \infty$. Since 
\[\norminf{wp_n}=\max\{\lambda_k:N_n\leq k <N_{n+1}\},\]
we have that $\norminf{wp_n}\rightarrow 0$ as $n \rightarrow \infty$. 

Conversely, suppose there exists a projection $\tilde{c}$ with the desired properties and $\norminf{wp_n}\rightarrow 0$. Let $w_N:=w\summ{n=1}{N}p_n$. Then $M_{w_N}$ is a finite rank operator for each $N\in \mathbb{N}$. Furthermore, using \cite[Proposition 5]{key-Han15} and the fact that the $p_n$'s are mutually orthogonal central projections, we obtain
\begin{eqnarray*}
\norm{M_w-M_{w_N}}=\norminf{w-w_N} =\norminf{w\summ{n=N+1}{\infty}p_n}=\supu{n>N}\norminf{wp_n}\rightarrow 0.
\end{eqnarray*}
It follows that $M_w$ is the limit of a sequence of finite rank operators and hence compact.
\end{proof}
\end{theorem}

Next, we consider multiplication operators between Orlicz spaces.

\begin{theorem}\label{T1 23/08/17}
Suppose $\varphi_i$ ($i=1,2,3$) are Orlicz functions with $\varphi_3\in \Delta_2$ and suppose one of the following conditions holds
\begin{enumerate}
\item either there exists some $k$ such that $\varphi^{-1}_1(t)\varphi_3^{-1}(t)k^{-1} \leq \varphi_2^{-1}(t)\leq k \varphi^{-1}_1(t)\varphi_3^{-1}(t)$, for all $t\geq 0$; 
\item or there exists an Orlicz function $\psi$ such that $\varphi_3=\psi^* \circ \varphi_2$ and $\psi \circ \varphi_2=\varphi_1$, with $\varphi_2$ satisfying the $\nabla'$-condition.
\end{enumerate}
If $w\in S(\mathcal{A},\tau)$, then $M_w$ is a compact multiplication operator from $L^{\varphi_1}(\tau)$ into $L^{\varphi_2}(\tau)$ if and only if $w\in L^{\varphi_3}(\tau)$ and there exists a sequence $(p_n)_{n=1}^\infty$ of mutually orthogonal central projections such that $w\tilde{c}=w$, where $\tilde{c}=\summ{n=1}{\infty}p_n$ and $\tilde{c}\mathcal{A}=\dsumx{n=1}{\infty} \mathcal{A}_n$, where each $p_n\mathcal{A}=\mathcal{A}_n$ is a finite type 1 factor.
\begin{proof}
Let $w\in L^{\varphi_3}(\tau)$ and suppose there exists a projection $\tilde{c}$ with the requisite properties. Let $w_N=w\summ{n=1}{N}p_n$. Then $M_{w_N}$ is a finite rank operator for each $N\in \mathbb{N}$ and $w-w_N\in L^{\varphi_3}(\tau)$, since $w_N\in L^1\cap L^{\infty}(\tau)\subseteq L^{\varphi_3}(\tau)$ for each $N\in \mathbb{N}$. It follows from either of Theorem \ref{T1 07/09/17} or Theorem \ref{T1 25/08/17} (depending on whether condition (1) or (2) holds) that $M_w-M_{w_N}=M_{w-w_N}$ is a bounded multiplication operator from $L^{\varphi_1}(\tau)$ into $L^{\varphi_2}(\tau)$ with 
\begin{eqnarray}
\norm{M_w-M_{w_N}}\leq 2\norm{w-w_N}_{\varphi_3}. \label{e2 23/08/17}
\end{eqnarray}
Note that $\summ{n=N+1}{\infty}p_n=\tilde{c}-\summ{n=1}{N}p_n\searrow 0$. So, on using the fact that each $p_n$ is a central projection (and hence also $\summ{n=N+1}{\infty}p_n$), we obtain 
\[w-w_N=w\summ{n=N+1}{\infty}p_n=w^{1/2}\left(\summ{n=N+1}{\infty}p_n \right)w^{1/2}\searrow 0.\]
Since $\varphi_3\in \Delta_2$, $L^{\varphi_3}(\tau)$ has absolutely continuous norm and therefore $\norm{w-w_N}_{\varphi_3}\searrow 0$. Using (\ref{e2 23/08/17}) this implies that $M_w$ is the limit of a sequence of finite rank operators and hence compact.

Conversely, if $M_w$ is a compact operator, then by Theorem \ref{T4.2 Lab14}, there exists a central projection $\tilde{c}$ with the desired properties. Furthermore, since $M_w$ is a compact operator, it is bounded and therefore $w\in L^{\varphi_3}(\tau)$, by either Theorem \ref{T1 07/09/17} or Theorem \ref{T1 25/08/17} (depending on whether condition (1) or (2) holds).
\end{proof}
\end{theorem}

Suppose $1< q < p$ with $r$ such that $1/p+1/r=1/q$. Then for $\varphi_1(t):=t^p$, $\varphi_2(t)=t^q$ and $\varphi_3(t)=t^r$, we trivially have \[\varphi^{-1}_1(t)\varphi_3^{-1}(t)= \varphi_2^{-1}(t), \qquad \forall  t\geq 0.\] 
We therefore obtain the following corollary.

\begin{corollary}
Let $1< q<p<\infty$ and $w\in S(\mathcal{A},\tau)^+$. Then $M_w$ is a compact multiplication operator from $L^p(\tau)$ into $L^q(\tau)$ if and only if 
\begin{enumerate}
\item there exists a sequence $(p_n)_{n=1}^\infty$ of mutually orthogonal central projections such that $w\tilde{c}=w$, where $\tilde{c}=\summ{n=1}{\infty}p_n$, $\tilde{c}\mathcal{A}=\dsumx{n=1}{\infty} \mathcal{A}_n$ and each $\mathcal{A}_n=p_n\mathcal{A}$ is a finite type 1 factor;
\item and $w\in L^{r}(\tau)$, where $1/p+1/r=1/q$.
\end{enumerate}
\end{corollary}

We finish by characterizing the compactness of multiplication operators from $L^p(\tau)$ into $L^q(\tau)$ for the case $1\leq p<q<\infty$.

\begin{theorem} \label{T1b 03/08/17}
Let $\mathcal{A}\subseteq \mathcal{B}(H)$ be a semi-finite von Neumann algebra equipped with a faithful normal semi-finite trace $\tau$, $1\leq p<q<\infty$ and $w\in S(\mathcal{A},\tau)^+$. Then $M_w$ is a compact multiplication operator from $L^p(\tau)$ into $L^q(\tau)$ if and only if 
\begin{enumerate}
\item there exists a sequence $(p_n)_{n=1}^\infty$ of mutually orthogonal central projections such that $w\tilde{c}=w$, where $\tilde{c}=\summ{n=1}{\infty}p_n$, $\tilde{c}\mathcal{A}=\dsumx{n=1}{\infty} \mathcal{A}_n$ and each $\mathcal{A}_n=p_n\mathcal{A}\cong B(H_n)$ with $H_n$ finite-dimensional;
\item and $\frac{\norminf{wp_n}}{k_n^{1/s}}\rightarrow 0$, where $1/q+1/s=1/p$ and $k_n$ is the value of the trace of a minimal projection in $B(H_n)$.
\end{enumerate}
\begin{proof}
Suppose conditions $(1)$ and $(2)$ hold. Then $M_w$ is a bounded multiplication operator from $L^p(\tau)$ into $L^q(\tau)$, by Theorem \ref{T3 26/10/17}. If we let $w_N:=w\summ{n=1}{N}p_n$, then for each $N\in \mathbb{N}$, $M_{w_N}$ is a finite rank operator. Furthermore, using Theorem \ref{T3 26/10/17}, we obtain \[\norm{M_w-M_{w_N}}=\norm{M_{w\summ{n=N+1}{\infty}p_n}}=\supu{n>N}\frac{\norminf{wp_n}}{k_n^{1/s}}\rightarrow 0.\] The operator $M_w$ is therefore the limit of a sequence of finite rank operators and hence compact. 

Conversely, if $M_w$ is compact, then, by Theorem \ref{T4.2 Lab14}, there exists a central projection $\tilde{c}$ with the desired properties. 
Since $\mathcal{A}_n$ is a finite type 1 factor, we can write  $wp_n=\summ{k=N_{n-1}+1}{N_n}\lambda_kq_k$, where the $\lambda_k$ are eigenvalues 
of $wp_n$ (repeated according to multiplicity) and the $q_k$ are projections onto $1$-dimensional subspaces of $p_n(H)$ with $q_k q_m=0$ 
if $k\neq m$. Consider the sequence $\{x_k\}_{k=1}^\infty$, where $x_k:=\tau(q_k)^{-1/p}q_k \in L^p(\tau)$. Recall that 
$L^p(\tau)^*\cong L^{p'}(\tau)$, where $1/p+1/p'=1$, and that the isometric isomorphism is given by $y \mapsto \tau(y\cdot)$.  If 
$y\in L^{p'}(\tau)$, then on using Lemma \ref{L2 20/10/17} at appropriate points, we have that
\begin{eqnarray}
|\ip{y}{x_k}|&=&|\tau(q_k)^{-1/p}\tau(yq_k)| \nonumber \\
&\leq&\tau(q_k)^{-1/p}\normx{yq_k}{1} \nonumber \\
&=&\tau(q_k)^{-1/p}\norminf{yq_k}\tau(q_k) \qquad \nonumber \\
&=&\norminf{yq_k}\tau(q_k)^{1/p'} \nonumber \\
&=&\normx{yq_k}{p'}  \qquad \label{e2 T1b 03/08/17}
\end{eqnarray}
Note that, by Lemma \ref{LP1 10/10/17}, $\normx{y}{p'}^{p'}=\summ{n=1}{\infty}\tau(|yp_n|^{p'})$. It follows that $\normx{yp_n}{p'}^{p'}=\tau(|yp_n|^{p'})\rightarrow 0$. Since $\normx{yq_k}{p'}\leq \normx{yp_n}{p'}$ for $N_{n-1}<k\leq N_n$, this implies that $\normx{yq_k}{p'}\rightarrow 0$. In view of the fact that $y\in L^{p'}(\tau)$ was arbitrary, this will, by means of (\ref{e2 T1b 03/08/17}), imply that $x_k\rightarrow 0$ weakly. For $N_{n-1}<k\leq N_n$, we have
\begin{eqnarray}
\frac{\lambda_k}{k_n^{1/s}}&=&\frac{\norminf{wq_k}}{k_n^{1/s}}=\frac{\normx{wq_k}{q}\tau(q_k)^{-1/q}}{\tau(q_k)^{1/p-1/q}}\nonumber \\
&=&\normx{wq_k \tau(q_k)^{-1/p}}{q}=\normx{M_w x_k}{q}. \label{e1 T1b 03/08/17}
\end{eqnarray}
Given that compact operators map weakly convergent sequences onto norm convergent sequences, the weak convergence of $x_k$ to $0$ ensures that $\normx{M_w x_k}{q}\rightarrow 0$.  Since 
\begin{eqnarray*}
\frac{\norminf{wp_n}}{k_n^{1/s}}&=&\frac{\max\{\lambda:\text{$\lambda$ is an eigenvalue of $wp_n$}\}}{k_n^{1/s}} \\
&=&\frac{\max\{\lambda_k:N_{n-1}<k\leq N_n\}}{k_n^{1/s}} \\
&=&\max\{ \normx{M_w x_k}{q}:N_{n-1}<k\leq N_n\} \qquad 
\end{eqnarray*}
(where we used (\ref{e1 T1b 03/08/17}) to obtain the last equality), it follows that $\frac{\norminf{wp_n}}{k_n^{1/s}}\rightarrow 0$. 
\end{proof}
\end{theorem}

\section{Acknowledgements}

The contribution of L. E. Labuschagne is based on research partially supported by the National Research Foundation (IPRR Grant 96128). Any opinion, findings and conclusions or recommendations expressed in this material, are those of the author, and therefore the NRF do not accept any liability in regard thereto.

\section*{References}

\bibliographystyle{amsplain}

\begin{thebibliography}{00}
\bibitem{key-Alfsen03} E.M. Alfsen and F.W. Shultz, \emph{Geometry of state spaces of operator algebras}, Mathematics: Theory and applications, Springer Science and Business Media, 2003.
\bibitem{key-Bala13} P. Bala, A. Gupta and N. Bhatia, {Multiplication operators on Orlicz-Lorentz sequence spaces}, \emph{Int. J. Math. Anal.}, \textbf{7}(30)(2013), 1461 - 1469.
\bibitem{key-BS88} C. Bennett and R. Sharpley, \emph{Interpolation of operators}, Academic Press, 1988.
\bibitem{key-Black06} B. Blackadar, \emph{Operator algebras: Theory of $C^*$-algebras and von Neumann algebras}, Encyclopaedia of Mathematical Sciences, Volume 122, Operator algebras and non-commutative geometry $III$, Springer-Verlag, 2006.
\bibitem{key-Cala08}J.M. Calabuig, O. Delgado and E.A. S\'{a}nchez P\'{e}rez,  Generalized perfect spaces, \emph{Indag. Math.}, \textbf{19} (2008), 359-378.
\bibitem{key-Chaw16} T. Chawziuk, Y. Cui, Y. Estaremi, H. Hudzik and R. Kaczmarek, Composition and multiplication operators between Orlicz function spaces, \emph{Journal of Inequalities and Applications}(2016): Paper No. 52, 18 pp.
\bibitem{key-Chaw17} T. Chawziuk, Y. Estaremi, H. Hudzik, S. Maghsoudi and I. Rahmani,  {Basic properties of multiplication and composition operators between distinct Orlicz spaces}, \emph{Rev. Mat. Comput.}, \textbf{30}(2017), 335-367.
\bibitem{key-Pag} B. de Pagter, {Non-commutative Banach function spaces}, \emph{Positivity},  197-227. In: \emph{Trends Math.}, Birkh\"auser, Basel, 2007.
\bibitem{key-Delgado10} O. Delgado and E.A. S\'{a}nchez P\'{e}rez, Summability properties for multiplication operators on Banach function spaces, \emph{Integr. Equ. Oper. Theory}, \textbf{66} (2010), 197-214.
\bibitem{key-Diestel95} J. Diestel, H. Jarchow and A. Tonge, \emph{Absolutely summing operators}, Cambridge University Press, 1995.
\bibitem{key-dP89}P.G. Dodds, T.K.-Y. Dodds, and B. de Pagter, Non-commutative Banach function spaces, \emph{Math. Z.}, \textbf{201}(1989), 583-597.
\bibitem{key-DDP2} P.G. Dodds, T.K.-Y. Dodds and B. de Pagter, {Non-commutative K\"othe duality}, \emph{Trans. Amer. Math. Soc.}, \textbf{339}(1993), 717-750.
\bibitem{key-Dodds14} P.G. Dodds and B. de Pagter, Normed K\"othe spaces: A non-commutative viewpoint, \emph{Indag. Math.}, \textbf{25}(2014), 206-249.
\bibitem{key-DDS14} P.G. Dodds, T.K.-Y. Dodds and F.A. Sukochev, On $p$-convexity and $q$-concavity in non-commutative symmetric spaces, \emph{Integr. Equ. Oper. Theory}, \textbf{78} (2014), 91-114.
\bibitem{key-Fack86} T. Fack and H. Kosaki, {Generalized $s$-numbers of $\tau$-measurable operators}, \emph{Pacific J. Math.}, \textbf{123}(1986), 269-300.
\bibitem{key-Han15}Y. Han, Products of noncommutative Calder\'{o}n-Lozanovski\u{i},  \emph{Math. Inequal. Appl.} \textbf{18} (2015), no. 4, 1341–1366
\bibitem{key-Han16}Y. Han, Generalized duality and product of some noncommutative symmetric spaces, \emph{Internat. J. Math.}, \textbf{27} (2016), no.9, 21 pages.
\bibitem{key-Hoa14} D. T. Hoa, Some inequalities for measurable operators, \emph{Int. J. Math. Anal.}, \textbf{8}(2014)(22), 1083-1087.
\bibitem{key-K1} R.V. Kadison and J.R. Ringrose, \emph{Fundamentals of the theory of operator algebras, Volume 1}, Birkh\"auser, Academic Press, 1983.
\bibitem{key-K2} R.V. Kadison and J.R. Ringrose, \emph{Fundamentals of the theory of operator algebras, Volume 2, Advanced theory}, Birkh\"auser, Academic Press, 1983.
\bibitem{key-Kras61} M.A. Krasnosel'skii and Y.B. Rutickii, \emph{Convex functions and Orlicz spaces}, Noordhoff Ltd., 1961.
\bibitem{key-Lab13} L.E. Labuschagne and W.A. Majewski, Maps on noncommutative Orlicz spaces, \emph{Illinois J. Math.}, \textbf{55}(3)(2011), 1053-1081.
\bibitem{key-Lab14} L.E. Labuschagne, {Multipliers on noncommutative Orlicz spaces}, \emph{Quaest. Math.}, \textbf{37}(2014), 531-546.
\bibitem{key-Mali89} L. Maligranda and L.E. Persson, Generalized duality of some Banach function spaces, \emph{Indag. Math.}, \textbf{92} (1989), 323-338.
\bibitem{key-Mali10} L. Maligranda and E. Nakai, Pointwise multipliers of Orlicz spaces, \emph{Arch. Math.} \textbf{95} (2010), 251-256.
\bibitem{key-Ramos17} J.C. Ramos-Fern\'{a}ndez and M. Salas-Brown, {On multiplication operators acting on K\"othe sequence spaces}, \emph{Afr. Mat.}, \textbf{28} (2017), 661-667.
\bibitem{key-Rao91} M.M. Rao and Z.D. Ren, \emph{Theory of Orlicz spaces}, Marcel Dekker Inc., 1991.
\bibitem{key-Suk16} F.A. Sukochev, H\"older inequality for symmetric operator spaces and trace property of $K$-cycles, \emph{Bull. Lon. Mat. Soc.}, \textbf{48}(4)(2016), 637-647. 
\bibitem{key-Takagi99} H. Takagi and K. Yokouchi, Multiplication and composition operators between two $L^{p}$- spaces,  American Mathematical Society, \emph{Contemp. Math.}, \textbf{232}(1999), 321-338.
\bibitem{key-Tak79}M. Takesaki, \emph{Theory of operator algebras I and II}, Springer-Verlag, New York, 1979.
\bibitem{key-Terp1}M. Terp, \emph{$L^p$-spaces associated with von Neumann algebras}, Rapport No. 3a, University of Copenhagen, (1981).
\end{thebibliography}

\end{document}